\documentclass[10pt]{article}

\usepackage{amsmath,amsthm,amssymb,latexsym}

\begin{document}


\newcommand{\holdreg}{\sigma}
\newcommand{\newholdreg}{\gamma}
\newcommand{\reg}{m}
\newcommand{\newreg}{k}
\newcommand{\basereg}{r}


\newcommand{\eg}{g_+}               
\newcommand{\eRic}{Ric_+}           
\newcommand{\eS}{S_+}               
\newcommand{\eLap}{\Delta_+}        
\newcommand{\deff}{\hat{\rho}}      
\newcommand{\cg}{\hat{g}}           


\newcommand{\geodeff}{\tilde{\rho}}     
\newcommand{\geog}{\tilde{g}}           
\newcommand{\geocoord}{\tilde{x}}       
\newcommand{\geocoder}{\tilde{\nabla}}  
\newcommand{\geoLap}{\tilde{\Delta}}    
\newcommand{\geoN}{\tilde{\mathcal{N}}} 
\newcommand{\geoL}{\tilde{L}}           
\newcommand{\geoRic}{\tilde{R}ic}       
\newcommand{\geoS}{\tilde{S}}           
\newcommand{\geoA}{\tilde{A}}           
\newcommand{\geoChr}{\tilde{\Gamma}}    
\newcommand{\geoinv}{\tilde{g}^{-1}}    
\newcommand{\georoot}{(\tilde{g}^{00})^{-\frac{1}{2}}}


\newcommand{\csbackdeff}{\rho}      
\newcommand{\csbackg}{g}            
\newcommand{\csbackRic}{Ric}        
\newcommand{\csbackS}{S}            
\newcommand{\csbackLap}{\Delta}     
\newcommand{\csbackD}{\mathcal{L}}  
\newcommand{\csbackcoder}{\nabla}   
\newcommand{\csbackChr}{\Gamma}     


\newcommand{\csdeff}{\rho}          
\newcommand{\csg}{g}                
\newcommand{\csRic}{Ric}            
\newcommand{\csS}{S}                
\newcommand{\csLap}{\Delta}         
\newcommand{\cscoder}{\nabla}       
\newcommand{\csChr}{\Gamma}         
\newcommand{\csA}{A}                
\newcommand{\csinv}{g^{-1}}         
\newcommand{\csroot}{(g^{00})^{-\frac{1}{2}}}


\newcommand{\ord}{\mathrm{ord}}     
\newcommand{\grad}{\mathrm{grad}}   
\newcommand{\adj}{\mathrm{adj}}     
\newcommand{\tr}{\mathrm{tr}}       
\newcommand{\tf}{\mathrm{tf}}       
\newcommand{\N}{\mathcal{N}}        
\newcommand{\E}{\mathcal{E}}        
\newcommand{\D}{\mathcal{L}}        
\newcommand{\F}{\mathcal{F}}        
\newcommand{\T}{\mathcal{T}}        
\newcommand{\A}{\mathcal{A}}	
\newcommand{\B}{\mathcal{B}}	
\renewcommand{\P}{\mathcal{P}}        

\newtheorem{main}{Theorem}
\renewcommand{\themain}{\Alph{main}}
\newtheorem{thm}{Theorem}[section]
\newtheorem{lemma}[thm]{Lemma}
\newtheorem{prop}[thm]{Proposition}



\title{Boundary Regularity for Conformally Compact Einstein Metrics in
  Even Dimensions}
\author{Dylan William Helliwell}
\date{\today}

\maketitle


\abstract{
We study boundary regularity for conformally compact Einstein metrics in even dimensions by
generalizing the ideas of Michael Anderson found in \cite{Anderson} and \cite{Anderson2}.
Our method of approach is to view the vanishing of the Ambient
Obstruction tensor as an $n$th order system of equations for the
components of a compactification of the given metric.  This, together
with boundary conditions that the compactification is shown to satisfy
provide enough information to apply classical boundary regularity
results.  These results then provide local and global versions of
finite boundary regularity for the components of the compactification.
}


\section{Introduction}

In recent years, both the mathematics and theoretical physics
communities have shown a great deal of interest in studying the
analysis and geometry of conformally compact Einstein metrics.  In
particular, much progress has been made recently with regard to
boundary regularity of these metrics.  This paper provides one
approach for investigating this topic and studies the problem of
regularity in the even dimensional case.

The concept of a conformally compact metric was developed in 1963
by Roger Penrose \cite{PR}.  Since then, such metrics
have been studied in \cite{LeBrun, MazMel, Mazthesis, GLee}
along with many others; see
\cite{Lee} and the references therein, for example.  The physics
community has also become interested in conformally compact Einstein
metrics since the introduction of the so called AdS/CFT
correspondence by Maldacena \cite{Mal}.  See for example
\cite{Witten, HS, dHSS} along with the references therein.

An example stemming from the developments in physics, which shows how
the issue of boundary regularity of a conformally compact Einstein
metric arises, is the volume renormalization of a conformally compact
Einstein manifold.  This renormalization gives
rise to an invariant associated with the conformal boundary manifold
\cite{green}. In order to obtain the result, there must be some
compactification for the Einstein metric which is smooth enough at
the boundary for certain calculations to be made.

The first result concerning boundary regularity of conformally compact
Einstein metrics was negative.  Namely, Fefferman and Graham
\cite{FeffG} showed that if $M$ is odd dimensional, there are examples
where there is no smooth compactification despite smoothness of the
boundary metric.  More recently, though, positive results have been
proved.  In \cite{Anderson}, Anderson studies the four dimensional
case by considering the Bach tensor of a compactification $g = \rho^2
\eg$ with constant scalar curvature.  The Bach tensor is a classically
known natural tensor depending on two derivatives of curvature which
is conformally invariant in dimension four, and vanishes for Einstein
metrics.  Making use of these facts and working in special harmonic
coordinates for $g$, Anderson generates second order uniformly
elliptic equations for components of the Ricci tensor of $g$, and for
the components of $g$ itself.  He also derives boundary equations for
this system.  A bootstrap argument is then applied to prove boundary
regularity.  In \cite{Anderson2}, he presents a revised version of the
argument based on viewing the Bach equation, combined with his
boundary equations, as a fourth order elliptic boundary value problem
in the sense of \cite{ADN2, Morrey} for the components of the constant
scalar curvature compactification.

In another recent paper, Chru\'{s}ciel, Delay, Lee, and Skinner \cite{cdls} study boundary
regularity in general dimension by applying the uniformly degenerate
theory of \cite{AndChr} directly to a gauge broken Einstein equation.
They show that in the case when the boundary metric is $C^{\infty}$,
there is a compactification which, in suitable coordinates, is
$C^{\infty}$ up to the boundary in even dimensions, and which has an
asymptotic expansion involving logarithms in odd dimensions.

In this paper, we follow Anderson's first approach to prove finite boundary
regularity in general even dimensions.  The Bach tensor is not
conformally invariant in dimensions other than four, but there is a
generalization of the Bach tensor in each higher even dimension called
the ambient obstruction tensor \cite{FeffG, GHir}.  Like the Bach
tensor, it is conformally invariant and vanishes for Einstein
metrics.  We choose a constant scalar curvature compactification for
the given Einstein metric and we work in harmonic coordinates.  In this
setting, the vanishing of the ambient obstruction tensor will provide
us with a system of equations that are uniformly elliptic up to the
boundary.  However, these equations will be $n$th order and this
brings up the difficulty of finding the necessary boundary equations
to make a well posed boundary value problem.  The natural boundary data
is that of a boundary metric, and nothing more.  A large task, then,
is to derive new boundary equations that the compactification
satisfies and that make a well posed problem.  Once such a boundary
value problem is derived, classical theorems may be applied to yield
boundary regularity results.

We focus on metrics that satisfy the specific Einstein equation $\eRic
= -(n-1) \eg$ and henceforth, unless mentioned otherwise, this is
always what we mean when we say a metric is Einstein.  We say that
$\eg$ is $C^{\reg,\holdreg}$ conformally compact if for any
$C^{\reg+1,\holdreg}$ defining function, the resulting compactified
metric $g = \rho^2 \eg$ is $C^{\reg,\holdreg}$ up to the boundary
(which makes sense as long as $\overline{M}$ has a
$C^{\reg+1,\holdreg}$ structure).  The local and global versions of
our main result are as follows:

\begin{main} \label{mainlocalthm} \label{MAINLOCALTHM}
Let $\overline{M}$ be an $n$-dimensional $C^{\infty}$ manifold with
boundary, $n \geq 4$ and even.  Let $p \in \partial M$ and let $U
\subset \overline{M}$ be a neighborhood of $p$ with boundary portion
$D = U \cap \partial M$.  For $\basereg \geq n$ and $0 < \holdreg < 1$, let $\eg$ be a
$C^{\basereg-1,\holdreg}$ conformally compact Einstein metric on $U \cap M$.
Suppose that the conformal infinity of $\eg$ contains a metric $h \in
C^{\newreg,\newholdreg}(D)$, where $\newreg \geq \basereg$ and $0 <
\newholdreg < 1$.  Given a $C^{\infty}$ coordinate system on the
boundary in a neighborhood of $p$, there is a coordinate system on a
neighborhood $V \subset U$ of $p$ that is $C^{\basereg, \holdreg}$ compatible
with the given smooth structure and which restricts to the given
coordinates on the boundary, and there is a defining function $\rho
\in C^{\basereg-1,\holdreg}(V)$ in the new coordinates, such that $\rho^2
\eg$ has boundary metric $h$ and in the new coordinates, $\rho^2 \eg
\in C^{\newreg,\newholdreg}(V)$.
\end{main}

\begin{main} \label{mainglobalthm} \label{MAINGLOBALTHM}
Let $\overline{M}$ be a compact $n$-dimensional $C^{\infty}$ manifold
with boundary, $n \geq 4$ and even.  For $\basereg \geq n$ and $0 < \holdreg < 1$,
let $\eg$ be a $C^{\basereg-1,\holdreg}$ conformally compact Einstein metric on
$M$.  Suppose that the conformal infinity of $\eg$ contains a metric
$h \in C^{\newreg,\newholdreg}(\partial M)$, where $\newreg \geq \basereg$
and $0 < \newholdreg < 1$.  Then there is a
 $C^{\basereg, \holdreg}$ diffeomorphism
$\Psi: \overline{M} \longrightarrow \overline{M}$ which restricts to the identity on the boundary and there is a defining function $\check{\rho} \in C^{\basereg-1,\holdreg}(\overline{M})$ such that
$\check{\rho}^2\Psi^*(\eg)$ has boundary metric $h$ and
$\check{\rho}^2 \Psi^*(\eg) \in C^{\newreg,\newholdreg}(\overline{M})$.
\end{main}

We note that the coordinates in Theorem \ref{mainlocalthm}
depend on $\eg$ and $h$, but not their regularity.  Hence, it follows 
that if $h$ is $C^{\infty}$, then the new coordinates are $C^{\infty}$, and $\rho^2 \eg$ is $C^{\infty}$ in these new coordinates, for some defining function $\rho \in
C^{\basereg-1,\holdreg}(\overline{M})$.  Similarly, if $h$ is $C^{\infty}$ in Theorem \ref{mainglobalthm}, then $\check{\rho}^2 \Psi^*(\eg)$ is $C^{\infty}$ for some
$\check{\rho} \in C^{\basereg-1,\holdreg}(\overline{M})$.

Anderson argues in \cite{Anderson} and \cite{Anderson2} that a version
of these results holds for $n = 4$ when a conformal compactification
of $\eg$ is $L^{2,p}$ for some $p > 4$, instead of $C^{\basereg-1,\holdreg}$, $r \geq 4$.
I am unable to verify the theorem with this hypothesis.  Also,
Anderson's statement does not mention the change of coordinates, nor
the fact that defining function may not be smooth in the new
coordinates.

Observe that in Theorem \ref{mainlocalthm}
we cannot conclude that $\eg$ is $C^{\newreg,\newholdreg}$ conformally
compact in the new coordinates, since the defining function need not
be $C^{\newreg+1,\newholdreg}$ in the new coordinates.  Similarly, we cannot conclude that
$\Psi^*(\eg)$ is $C^{\newreg, \newholdreg}$ conformally compact in
Theorem \ref{mainglobalthm}.  This is a
consequence of our specific compactification as well as the change of
coordinates.  I expect that $\eg$ is $C^{\newreg,\newholdreg}$
conformally compact in the new coordinates.  This can be reduced to a
regularity question for the singular Yamabe problem, since as an
Einstein metric, $\eg$ is a (singular) constant scalar curvature
metric in its conformal class.  Using results of \cite{AndChr} on this
problem, we are able to conclude the following theorem:

\begin{main}\label{maincptthm} \label{MAINCPTTHM}
Let $\overline{M}$ be a compact $n$-dimensional $C^{\infty}$ manifold
with boundary, $n \geq 4$ and even.  For $\basereg \geq n$ and $0 < \holdreg < 1$, let
$\eg$ be a $C^{\basereg-1,\holdreg}$ conformally compact Einstein metric on
$M$.  Suppose that the conformal infinity of $\eg$ contains a metric
$h \in C^{\newreg,\newholdreg}(\partial M)$, where $\newreg \geq \basereg$ and $\newreg
\geq n+1$, and $0 < \newholdreg < 1$.  Then there is a $C^{\basereg, \holdreg}$ diffeomorphism
$\Psi: \overline{M} \longrightarrow \overline{M}$ which restricts to the identity on
$\partial M$ such that $\Psi^*(\eg)$ is $C^{\newreg,\newholdreg'}$ conformally compact for
some $\newholdreg'$, $0 < \newholdreg' \leq \newholdreg$.
\end{main}

Similarly to the situation with Theorems \ref{mainlocalthm} and
\ref{mainglobalthm}, this result allows us to conclude that $\Psi^*(\eg)$ is
$C^{\infty}$ conformally compact in the new smooth structure if $h$ is
$C^{\infty}$.

Theorem \ref{maincptthm} is deficient in that, when $\basereg = n$, it should also hold with
$\newreg = n$ and in general, it should hold with  $\newholdreg' = \newholdreg$, and there should be an analogous local result.  These deficiencies are a direct consequence
of corresponding deficiencies in the regularity theorem of
\cite{AndChr}.  Hence, improvement of their theorem would result in
improvement of Theorem \ref{maincptthm}.

The regularity theorem of \cite{cdls} using uniformly degenerate
methods assumes that $\eg$ is defined in a collar neighborhood of the full
boundary, $\eg$ is $C^2$ conformally compact, and $h \in
C^{\infty}(\partial M)$.  Our use of high order elliptic equations and
boundary conditions has the disadvantage that it seems to necessitate
the high a priori regularity assumptions in Theorems
\ref{mainlocalthm}--\ref{maincptthm}.  On the other hand, Theorem
\ref{mainlocalthm} is local, and all three theorems apply for boundary
metrics of finite regularity.

The outline of this paper is as follows.
We begin, in Section \ref{Background}, by developing the mathematical
background for studying conformally compact Einstein metrics.  We
define precisely what a conformally compact metric is and introduce a
number of related notions.  We also discuss boundary adapted harmonic
coordinates and discuss their properties.  To finish
the section, we introduce the constant scalar curvature
compactification.
     
In Section \ref{Approximate}, we recall the notions of geodesic
compactification and geodesic coordinates.  The geodesic
compactification for a conformally compact Einstein metric in its
geodesic coordinates has a particularly simple form.  This then can be
used to express curvature tensors of the compactification at the
boundary purely in terms of the boundary metric.  Unfortunately,
working with this compactification in these coordinates introduces a
loss of regularity.  To avoid this loss, we construct an approximate
version of the geodesic compactification and use the approximate
geodesic coordinates of \cite{AndChr} for this compactification.  We then show that at the boundary, these ``almost
geodesic compactifications'' produce the same expressions as exact
geodesic compactifications for derivatives up to a fixed finite order
of curvature tensors in terms of the boundary metric.

Section \ref{BVPchap} is dedicated to the derivation of the boundary
value problem.  The equations on the interior will come from
expressing the ambient obstruction tensor for a constant scalar
curvature compactification in harmonic coordinates.  In deriving
boundary conditions, we will find that some of the equations come to
us naturally, while others require a careful analysis of curvature
tensors at the boundary.  The results for almost geodesic
compactifications will play a major role in this analysis.  Our first order boundary conditions are different from and simpler than those of Anderson.  To derive them, we will not need to require that our harmonic coordinates restrict to harmonic coordinates for $h$ on the boundary. 

Theorems \ref{mainlocalthm} and \ref{mainglobalthm} are proved in
Section \ref{regularitychap}.  To prove Theorem \ref{mainlocalthm}, we
treat the $n$th order boundary value problem as a sequence of second
order problems much like Anderson's approach in \cite{Anderson}.  An
analysis of the boundary equations shows that they are sufficient for
the application of boundary regularity results to these second order
equations.  Following this, an involved bootstrap yields the result.
With Theorem \ref{mainlocalthm} proved, Theorem \ref{mainglobalthm}
then follows by a patching argument and an approximation theorem.

Section \ref{defining} contains our discussion of regularity of the
defining function.  The regularity theorem of \cite{AndChr} asserts
the existence of an asymptotic expansion involving log terms for the
solution of the singular Yamabe problem.  We show that in the case
that the Yamabe metric is actually Einstein, then no log terms occur.
This along with Theorem \ref{mainglobalthm} is used to prove
Theorem \ref{maincptthm}.

Finally, section \ref{syschap} discusses an alternative approach to the
proof of Theorem \ref{mainlocalthm}, which is to treat the problem as
a single elliptic boundary value problem.  This is similar to Anderson's analysis in \cite{Anderson2}.

\vspace{.25in}

I would like to thank Robin Graham for his help on this paper.  The many discussions we had on these topics were indispensable and without his patience and guidance, this paper would not exist.


\section{Background} \label{Background}

To begin, we set up the situation in which we will be working primarily.
In doing so, principal definitions and notions are introduced, along
with a number of facts and tools that will
be used later.  We also introduce our notation and conventions that will be
followed throughout this paper.

\subsection{Preliminaries and Main Definitions}

For the most part, we will be working on or near the boundary of an
$n$-dimensional $C^{\infty}$ manifold with boundary $\overline{M}$,
with interior $M$, and $C^{\infty}$ boundary $\partial M$.  When
working locally, let $U$ be an open set in $\overline{M}$
containing a nonempty boundary portion $D = U \cap \partial M$.  In
general, the exact choice of $U$ is not crucial, and so as necessary,
and without explicit mention, $U$ may be shrunk to a smaller domain,
still containing part of the boundary, so that all relevant notions
are well defined.  A function $\rho$ on $U$ is a defining function for
$D$ if $\rho|_D = 0$, $d\rho|_D \neq 0$, and for convention, $\rho$ is
positive in $M$.

Let $\reg \in \mathbb{N} \cup \{0\}$, and $0 < \holdreg < 1$.  We say
a metric $\eg$ on $U \cap M$ is $C^{\reg,\holdreg}$ (resp. $C^\reg$,
$C^{\infty}$) \emph{conformally compact} if $\rho^2 \eg$ extends to a
$C^{\reg,\holdreg}$ (resp. $C^\reg$, $C^{\infty}$) metric on $U$,
where $\rho$ is a $C^{\reg+1,\holdreg}$ (resp. $C^{\reg+1}$,
$C^{\infty}$) defining function for $D$.  Note that this makes sense
as long as $\overline{M}$ has a $C^{\reg+1,\holdreg}$
(resp. $C^{\reg+1}$, $C^{\infty}$) structure.  For a given $\rho$, let
$g = \rho^2 \eg$ be the extended metric. Then we say $g$ is a
\emph{compactification} for $\eg$.  Letting $\iota:D \rightarrow U$ be
inclusion, we call $h = \iota^*g$ a \emph{boundary metric}.  The
equivalence class $[h]$ whose elements arise from various choices of
$\rho$ is called the \emph{conformal infinity} of $\eg$.

It follows from these definitions that for any conformally compact
metric $\eg$, the restriction to $\partial M$ of $|d \rho|_{\rho^2
\eg}$ is invariant with respect to the choice of $\rho$.  If this
invariant is constantly equal to 1 and $\reg \geq 2$, then the
sectional curvatures of $\eg$ all approach $-1$ as we approach the
boundary.  See \cite{Anderson} or \cite{Mazthesis} for details.  We
say that $\eg$ is \emph{asymptotically hyperbolic} on $U$ if $\eg$ is
conformally compact and $|d \rho|_{\rho^2 \eg}^2 = 1$ on $D$.  If
$\eg$ is $C^{\reg}$ conformally compact on $U$ and $\rho^{-1}(1 -
|d\rho|_{\rho^2 \eg}^2) \in C^{\reg}(U)$, then we say $\eg$ is
\emph{$C^{\reg}$ asymptotically hyperbolic} on $U$.  It is
straightforward to check that this definition is independent of the
choice of $C^{\reg+2}$ defining function.  Also observe that if $\eg$ is
asymptotically hyperbolic and $C^{\reg+1}$ conformally compact, then
it is $C^\reg$ asymptotically hyperbolic.

We study conformally compact Einstein metrics, by which we mean
conformally compact metrics $\eg$ satisfying the Einstein equation
$\eRic = -(n-1) \eg$.  Such metrics are asymptotically hyperbolic, but
we will sometimes need to use the stronger condition mentioned above.
An Einstein metric that is $C^{\reg,\holdreg}$ conformally compact will
be called a \emph{$C^{\reg,\holdreg}$ conformally compact Einstein
metric}.  Similarly, an Einstein metric which is
$C^{\reg}$ asymptotically hyperbolic will be called a \emph{$C^{\reg}$
asymptotically hyperbolic Einstein metric}.

Because we will be working with various metrics, objects associated to
a given metric will be adorned with the same symbol.  For example, the
Ricci curvature tensor associated to $\eg$ will be denoted $\eRic$.
If a symbol is unadorned, then it corresponds to the unadorned metric
$g$.  Also, we note here that we do not change notation to distinguish
between various representations of a metric in different coordinates.
As this has the potential to generate confusion, we make clear what is
happening whenever we are dealing with multiple coordinate systems at
once.

\subsection{Boundary Adapted Harmonic Coordinates} \label{harmonic}

A coordinate system is called a \emph{boundary adapted coordinate
system} if it has the property that one of the coordinate functions is
a defining function for the boundary.  Note that in such a
coordinate system, the remaining coordinates restrict to a coordinate
system on $\partial M$.

We will be working with such coordinate systems throughout,
so we introduce some conventions.  Let $\{x^{\alpha}\}$ be a boundary
adapted coordinate system.  Then, $x^0$ will always be a defining
function, while $x^i$ ($i \neq 0$) will always refer to the remaining
coordinates.  Consistent with this, any use of indices will follow the
convention that Roman indices will range from 1 to $n-1$, while Greek
indices range from 0 to $n-1$.

While working with tensors in coordinates and to represent terms of
order less than the principal part that appear in various equations,
we will often have a need to encapsulate terms depending on a certain
number of derivatives of the metric or other tensors.  With this in
mind, let $\P(\partial^{p_1} s_1, \partial^{p_2} s_2,\ldots)$ denote
an expression whose components are
polynomials in the components of the tensors $s_r$ and their coordinate
derivatives up to order $p_r$.  We will use ``$\partial_t$'' to represent
derivatives with respect to coordinates $x^i$, while
``$\partial_0$'' indicates differentiation only with respect to $x^0$.
Also, the use of ``$\nabla^{p} s$''
instead of ``$\partial^p s$'' will indicate
dependence on the components of covariant derivatives of $s$ up to order
$p$.  Finally, different instances of this notation will not mean the
same expression.

When there is no chance for confusion, ``$\partial_{\alpha}$'' or a
comma followed by ``$\alpha$'' will be used to denote coordinate
differentiation with respect to the coordinate function $x^{\alpha}$.
Covariant derivatives will be represented with a ``$\nabla_{\alpha}$''.

For later reference, we note that the Riemannian Laplacian of a
function $f$ with respect to a metric $g$ is given by $\Delta f =
\tr_g(\nabla^2 f)$ so that in coordinates, $\Delta f = g^{\alpha
\beta} \partial_{\alpha} \partial_{\beta} f + \P(g^{-1}, \partial g,
\partial f)$.  More generally, letting $\D_l f$ be the $l$-fold trace
of $2l$ coordinate derivatives of $f$, we have $\Delta^l f = \D_l f +
\P(g^{-1}, \partial^{2l-1} g, \partial^{2l-1} f)$.

A coordinate system on $\overline{M}$ is called a \emph{harmonic
coordinate system} for a metric $g$ if each of the coordinate
functions is harmonic with respect to $g$.  It will turn out that
boundary adapted harmonic coordinates will prove most beneficial when
it comes to questions of regularity.  There are two main reasons for
this.  First, the regularity of a metric is preserved when
transforming to its harmonic coordinates.  Second, the components of
the Ricci tensor can be expressed as $\D$ acting on components of the
metric plus lower order terms.  These facts were established and
existence of harmonic coordinates about a point was proved in
\cite{DetKaz}.  To account for a boundary, we have the following
proposition.

\begin{prop} \label{hcoordsprop}
Let $U \subset \overline{M}$ be a coordinate domain with boundary
portion $D$ and let $\{x^i\}$ be coordinates for $D$ which are compatible with the coordinates on $U$.  Let $g$ be a
metric in $C^{\reg,\holdreg}(U)$, where $\reg \geq 1$ and $0 <
\holdreg < 1$.  Then near any point $p \in D$, there exists a boundary
adapted harmonic coordinate system for $g$ that is
$C^{\reg+1,\holdreg}$ related to the given coordinates on $U$ and that restricts to
the given coordinates $\{x^i\}$ on the boundary.
\end{prop}

\begin{proof}
Let $\{y^{\alpha}\}$ be a $C^{\infty}$ boundary adapted coordinate
system in a neighborhood of $p$ such that when $y^0 = 0$, the
remaining coordinates $\{y^i\}$ restrict to $\{x^i\}$.  Consider an
open set $V$ with the following properties.  First, $V$ is the
interior of a $C^{\infty}$ manifold with boundary, diffeomorphic to an
open $n$-dimensional ball.  Second, $\partial V \cap \partial M$ is
diffeomorphic to a closed $(n-1)$-dimensional ball in $\partial M$, and
contains $p$ in its interior.  Third, $\overline{V}$ is in the domain of the
coordinate system $\{y^{\alpha}\}$.  Now, construct harmonic functions
$x^{\alpha}$ near $p$ by solving the Dirichlet problem
\[
\left\{ \begin{array}{l}
        \Delta x^{\alpha} = 0 \\
        x^{\alpha}|_{\partial V} = y^{\alpha}
        \end{array} \right.
\]
Existence of such functions can be found in \cite{G&T}.  Moreover, by the
maximum principle, we know $x^0 > 0$ in $V$, and by Hopf's Lemma, we
know that $d x^0|_{\partial V \cap \partial M} \neq 0$.  Hence, $x^0$
is a defining function for $\partial M$ near $p$.

The set of functions $\{x^0, x^i; i=1, \ldots, n-1\}$ will form
coordinates on an open set $W$ near $p$.  By elliptic boundary
regularity (see \cite{G&T}), these coordinate functions are
$C^{\reg+1,\holdreg}(W)$ with respect to $\{y^{\alpha}\}$, since $g$ is in
$C^{\reg,\holdreg}(W)$.
\end{proof}

Observe that $g$ is $C^{\reg,\holdreg}$ in these new coordinates since
viewing $g$ after a coordinate change involves only first derivatives
of the new coordinates.  Details related to this can be found in
\cite{DetKaz}.  Moreover, since the construction here does not alter
the boundary coordinates, the regularity of $h$ is preserved as well.

Besides preserving the regularity of $g$, harmonic coordinates also
produce a simplified formula for the Ricci tensor, which will turn out
to have a number of applications in this paper.

\begin{lemma}
In harmonic coordinates for $g$,
\begin{equation} \label{Ric}
Ric_{\alpha \beta} = -\frac{1}{2}\D g_{\alpha \beta}
		     + \P(g^{-1}, \partial g).
\end{equation}
\end{lemma}

This lemma and its proof can be found in \cite{DetKaz}.  It is
basically a computation and we indicate the key ingredients here.
First, we note that having a harmonic coordinate system with respect
to $g$ is equivalent to having traces of the Christoffel symbols all
equal to zero:
\[
\csChr^{\alpha} = \csg^{\beta \delta}
                        \csChr^{\alpha}_{\beta \delta}
		      = 0.
\]
We can differentiate this equation and use the resulting conditions
to reduce a standard formula for Ricci to \eqref{Ric}.

\subsection{Constant Scalar Curvature Compactification}

As indicated in the introduction, we will find that making an
appropriate choice of conformal compactification is crucial to
guaranteeing a maximal regularity result.  Here, we introduce the
constant scalar curvature compactification.

\begin{prop} \label{csprop}
Let $U \subset \overline{M}$ be an open set with boundary portion $D$
and let $\eg$ be $C^{\reg,\holdreg}$ conformally compact on $U \cap
M$, where $\reg \geq 2$ and $0 < \holdreg < 1$.  Let $h \in
C^{\reg,\holdreg}(D)$ be an element of the conformal infinity for
$\eg$.  Then near any point $p \in D$, there is a defining function
$\csbackdeff$ such that $\csbackg = \csbackdeff^2 \eg \in
C^{\reg,\holdreg}$ near $p$ with the properties that $\csbackg$ has
boundary metric $h$ and scalar curvature constantly equal to
$-n(n-1)$.  Moreover, if $\hat{\rho}$ is a $C^{\reg+1,\holdreg}$ defining
function such that $\hat{\rho}^2 \eg$ restricts to $h$ then $\csbackdeff/\hat{\rho} \in
C^{\reg,\holdreg}$ near $p$.
\end{prop}

\begin{proof}
We start by constructing a set similar to that used in the proof of
Proposition \ref{hcoordsprop}.  In particular, we consider an open set
$V \subset U \cap M$ with the following properties.  First,
$\overline{V} \subset U$.  Second, $V$ is the interior of a
$C^{\infty}$ manifold with boundary, diffeomorphic to an open
$n$-dimensional ball.  Third, $\partial V \cap \partial M$ is a subset
of $D$ diffeomorphic to a closed $(n-1)$-dimensional ball in $\partial
M$, and contains $p$ in its interior.  With this in place, let
$\hat{g} = \hat{\rho}^{2}\eg$ be a compactification with boundary
metric $h$, where $\hat{\rho} \in C^{\reg+1,\holdreg}(\overline{V})$ and
$\hat{g} \in C^{\reg,\holdreg}$.  Our goal is to find a positive
function $v$ on $\overline{V}$ such that $\csbackg = v^{\frac{4}{n-2}}
\hat{g}$ has constant scalar curvature.  Looking at how scalar
curvature changes under this conformal transformation, we have
\begin{equation} \label{confS}
\csbackS = \hat{S}\,v^{-\left(\frac{4}{n-2}\right)}
           + \frac{4(1-n)}{n-2}
             v^{-\left(\frac{n+2}{n-2}\right)}\hat{\Delta} v.
\end{equation}
Setting $\csbackS = -n(n-1)$ (any negative number will work) we may
reduce our problem to the following:
\[
\left\{ \begin{array}{l}
        \hat{\Delta} v - \frac{(n-2)\hat{S}}{4(n-1)} v
                 - \frac{n(n-2)}{4} v^{\frac{n+2}{n-2}} = 0 \\
        v|_{\partial V} = 1.
        \end{array} \right.
\]
This is a special case of the Yamabe problem on a manifold with
boundary and always has a $C^{2,\holdreg}$ solution \cite{Ma}.
Moreover, since the given differential operator is elliptic, for $\hat{g}
\in C^{\reg,\holdreg}$, a bootstrap procedure shows that any such
solution $v$ is in $C^{\reg,\holdreg}(\overline{V})$.

With this, we can construct a new defining function $\csbackdeff =
v^{\frac{2}{n-2}} \hat{\rho}$, with which to compactify.  Note that since
$v = 1$ on $D \cap \partial V$, the boundary metric $h$ is not
changed near $p$.
\end{proof}

We remark that Proposition \ref{csprop} has a global analogue.
Namely, if $\eg$ is $C^{\reg,\holdreg}$ conformally compact on $M$
then there is a constant scalar curvature compactification $\rho^2 \eg$ in
$C^{\reg,\holdreg}(\overline{M})$ and for any $C^{\reg+1,\holdreg}$ defining
function $\hat{\rho}$ for $\partial M$, we have
$\csbackdeff/\hat{\rho} \in C^{\reg,\holdreg}(\overline{M})$.  This
follows by using $M$ in place of $V$ in the proof above.

If $\eg$ is Einstein and we work in harmonic coordinates for the
constant scalar curvature compactification $\csbackg$, then interior
regularity for $\csbackg$ is significantly better than that given in
Proposition \ref{csprop}.  In fact, we have the following result.
Note that $\partial M$ is not involved, and $\csbackdeff$ need not be
a defining function.

\begin{prop} \label{intreg}
Let $V$ be an open subset of $M$ and let $\eg$ be an Einstein metric
in $C^{2}(V)$.  Let $\csbackg = \csbackdeff^2 \eg$ have constant
scalar curvature, where $\csbackdeff \in C^{2}(V)$ is a positive
function and let $W \subset V$ be an open set on which harmonic coordinates
for $\csbackg$ are defined.  Then in these coordinates, $\csbackg$ and
$\csbackdeff$ are both in $C^{\infty}(W)$.
\end{prop}

\begin{proof}
Writing out the equation $\eRic = -(n-1)\eg$ in terms of $\csbackg$,
we find
\[
\csbackRic_{\alpha \beta}
     = \P(\csbackg^{-1}, \partial \csbackg,
          \csbackdeff^{-1}, \partial^2 \csbackdeff).
\]
Since we are in harmonic coordinates, we can use \eqref{Ric} and
rearrange to get
\begin{equation} \label{Dgbar}
\csbackD \csbackg_{\alpha \beta}
     = \P(\csbackg^{-1}, \partial \csbackg,
          \csbackdeff^{-1}, \partial^2 \csbackdeff).
\end{equation}
Also, similar to \eqref{confS}, we have
\[
\eS = \csbackS \csbackdeff^2
      + 2(n-1)\csbackdeff \csbackLap \csbackdeff
      - n(n-1)|d\csbackdeff|_{\csbackg}^2.
\]
By isolating the second derivatives of $\csbackdeff$, this equation
reduces to
\begin{equation} \label {Drhobar}
\csbackD \csbackdeff = \P(\csbackg^{-1}, \partial \csbackg,
                          \csbackdeff^{-1}, \partial \csbackdeff)
\end{equation}
since the scalar curvatures of $\eg$ and $\csbackg$ are both constant.

Now we run a double bootstrap to get the result.  To start, note that
in harmonic coordinates for $\csbackg$, the function $\csbackdeff$ and the
components of $\csbackg$ are all in $C^{1,\holdreg}(W)$.  Proceeding
by induction, let $\reg \geq 1$ and suppose that $\csbackdeff$ and the
components of $\csbackg$ are all in $C^{\reg,\holdreg}(W)$.  Then the
right hand side of \eqref{Drhobar} is in $C^{\reg-1,\holdreg}(W)$, so by
elliptic regularity, we may conclude that $\csbackdeff \in
C^{\reg+1,\holdreg}(W)$.  Using this, we then observe that the right hand
side of \eqref{Dgbar} is in $C^{\reg-1,\holdreg}(W)$, and so we have
$\csbackg_{\alpha \beta} \in C^{\reg+1,\holdreg}(W)$.  Therefore, by
induction it follows that $\csbackg_{\alpha \beta}$ and $\csbackdeff$
are in $C^{\reg,\holdreg}(W)$ for all $\reg$.
\end{proof}


\section{Almost Geodesic Coordinates for an Almost Geodesic
  Compactification} \label{Approximate}

Geodesic coordinates provide a means of expressing the Taylor
expansion at the boundary of a geodesic
compactification of a conformally compact Einstein metric in terms of
the boundary metric.  Unfortunately, using exact geodesic coordinates
for an exact geodesic compactification results in a loss of
regularity.  Here, we construct approximate versions of the geodesic
compactification and geodesic coordinates which avoid the loss of
regularity.  The approximate geodesic coordinates that we use are the
approximate Gaussian coordinates introduced in \cite{AndChr} in the
case when the background metric is taken to be our ``almost geodesic
compactification''.  Also see \cite{cdls} for similar constructions.

\subsection{Geodesic Compactification and Geodesic Coordinates}

We start with a review of the geodesic compactification and its
associated geodesic coordinates.  This compactification and a number
of its properties were introduced in \cite{GLee}.

\subsubsection{Geodesic Compactification}\label{geocomp}

Given an asymptotically hyperbolic metric $\eg$ on a coordinate domain
$U \subset \overline{M}$ with boundary portion $D$, it is possible to
find a defining function $\geodeff$, with associated compactified
metric $\geog$, so that $|d \geodeff|_{\geog} = 1$ not just on $D$,
but on a neighborhood of the boundary.  Moreover, given a boundary
metric $h$ in the conformal infinity of $\eg$, $\geodeff$ can be
chosen so that the new compactification has $h$ as its boundary
metric.  To do this, let $g$ be a compactification for $\eg$, with
defining function $\rho$, such that its boundary metric is $h$.  To
determine $\geog$, write $\geodeff = e^u \rho$.  Then we want to find
$u$ such that $u|_{\partial M} = 0$ and $|d(e^u \rho)|^2_{e^{2u} g} =
1$ in $U$.  We may take boundary adapted coordinates $\{x^{\alpha}\}$
for $U$, and since $\rho$ is a $C^{\infty}$ defining function, we may
take it to be $x^0$.  In this setting, the constant length condition
gives:
\begin{equation} \label{const_length}
2g^{0 \alpha} u_{\alpha} + \rho |du|_g^2 = \frac{1-|d\rho|_g^2}{\rho}.
\end{equation}
We are guaranteed a solution to this equation by general theory of
first order partial differential equations.  As for regularity, if
$\eg$ is $C^{\reg}$ asymptotically hyperbolic, we know that the right
hand side of this equation is $C^{\reg}(U)$ and so the best we can say
for an exact solution $u$ to this differential equation is that $u \in
C^{\reg}(U)$.  Hence $\geodeff \in C^{\reg}(U)$ as well.  See
\cite{Lee2} for a detailed discussion of H\"{o}lder regularity results
for this equation.  We call this compactification a \emph{geodesic
compactification}, and the associated defining function is called a
\emph{geodesic defining function}.  The terminology here comes from
the fact that the integral curves of $\grad_{\geog}(\geodeff)$ are
geodesics.

\subsubsection{Geodesic Coordinates}

The geodesic defining function $\geodeff$ can be used to identify
points near $D$ with elements of $D \times [0,\varepsilon)$ and to split
$\geog$ as $\geog_{\geodeff}+d\geodeff^2$,
where $\geog_{\geodeff}$ is a one parameter family of metrics on $D$.
To achieve this, choose any coordinate system $\{\tilde{x}^i; i = 1,
\ldots, n-1\}$ on $D$.  Extend these coordinates inward by
keeping them constant on the integral curves of
$\grad_{\geog}(\geodeff)$.  These functions taken together with
$\geocoord^0 = \geodeff$ form a boundary adapted coordinate system.
Such coordinates are called \emph{geodesic coordinates} for $\geog$.

Regularity of $\geog$ is not preserved in its geodesic coordinates.
On the other hand, $\geog_{00} = 1$ and $\geog_{i0} = 0$ by
construction.  If $\geog$ is a geodesic compactification for an
Einstein metric $\eg$, then we can say more.  In particular, when $n
\geq 4$ is even and $\geog$ is sufficiently smooth with respect to its
geodesic coordinates, then for $0 \leq p \leq n-2$,
\begin{equation*}
(\partial_0^p \geog_{ij})|_D = \P(h^{-1}, \partial_t^p h).
\end{equation*}
Moreover, if $p$ is odd, then in fact $(\partial_0^p \geog_{ij})|_D =
0$.  It follows immediately that if we allow for derivatives with
respect to $\geocoord^i$, then
\begin{equation} \label{geocomps}
(\partial^p \geog_{ij})|_D = \P(h^{-1}, \partial_t^p h).
\end{equation}
These expressions, along with similar expressions when $n$ is odd, are
derived in \cite{green}.

\subsection{Almost Geodesic Compactification and Almost Geodesic
  Coordinates}

We now construct approximate versions of the geodesic
compactification and geodesic coordinates, which preserve the
useful expressions given above while avoiding the loss of regularity
inherent in exact geodesic coordinates.

\subsubsection{Error Terms} \label{error_notation}

Because we will be working with approximations, we introduce some
notation for the error terms that arise.  In a coordinate domain $U$,
by $\E_\reg$ we will denote a function which is in $C^\reg(U)$, and which is
$o(\rho^\reg)$ for some defining function $\rho$.  Note that a $C^\reg$
function $f$ is $\E_\reg$ if and only if all derivatives up to order $\reg$
are zero on the boundary.  Different instances of this notation will
not mean the same function.

\subsubsection{Almost Geodesic Defining Function and Almost Geodesic
  Compactification}

Our first step in producing an appropriate approximation to the
geodesic compactification is to construct the right defining function.
For a $C^{\reg}$ asymptotically hyperbolic metric $\eg$, an exact
solution to \eqref{const_length} is $C^{\reg}$.  If we allow an
approximate solution, we can do a bit better.  We say a function
$\geodeff$ is an \emph{almost geodesic defining function of order
$\reg$} for $D$ if $\geodeff$ is a defining function for $D$ in
$C^{\reg+1}(U)\cap C^{\infty}(U \cap M)$ and
$1-|d\geodeff|_{\geodeff^2 \eg}^2 = \geodeff\, \E_\reg$.  Construction
of a global almost geodesic defining function of order 2 can be found
in \cite{cdls}.  For higher orders and in the local setting, the
following lemma guarantees that such a function exists by using an
extension result in \cite{AndChr} in an inductive argument.

\begin{lemma} \label{ageodefflem}
Let $\eg$ be a $C^{\reg}$ asymptotically hyperbolic metric on a
coordinate domain $U$ with boundary portion $D$, and let $h \in
C^{\reg}(D)$ be an element of the conformal infinity for $\eg$.  Then
there exists an almost geodesic defining function $\geodeff$ of order
$\reg$ such that the boundary metric of $\geodeff^2 \eg$ is $h$.
\end{lemma}

\begin{proof}
Let $g = \rho^2 \eg$ be a $C^{\reg}$ compactification that restricts
to $h$ on the boundary, where $\rho$ is a $C^{\infty}$ defining
function.  Taking boundary adapted coordinates $\{x^{\alpha}\}$ with
$x^0 = \rho$, our goal is to find a function $u$ that solves
\eqref{const_length} approximately.  Let $f$ be the right hand side of
\eqref{const_length}, and note that since $\eg$ is $C^{\reg}$
asymptotically hyperbolic on $U$, $f \in C^{\reg}(U)$.  Formally
differentiating \eqref{const_length} $l-1$ times with respect to
$\rho$, $1 \leq l \leq \reg+1$ and setting $\rho$ equal to zero we
obtain
\[
(\partial_0^{l} u)|_D = \P_l
\]
where $\P_l = \P(g^{-1},(g^{00})^{-1}, \partial_0^{l-1} g,
\partial_0^{l-1} f, \partial_t \partial_0^{l-1} u)$.

Hence, starting with $v^0 = 0$, we can recursively define functions 
$v^l = \P_l \in C^{\reg-l+1}(D)$.  By Corollary 3.3.2 in
\cite{AndChr}, there is a function $v \in C^{\reg+1}(U) \cap
C^{\infty}(U \cap M)$ such that $(\partial_0^{l} v)|_D = v^l$.

Now, by construction $v$ solves \eqref{const_length} modulo $\E_\reg$, so
defining $\geodeff = e^v \rho$, the result follows.
\end{proof}

Given a $C^{\reg}$ asymptotically hyperbolic metric $\eg$, we say that a
metric $\geog$ is
an \emph{almost geodesic compactification of order $\reg$ associated to
$\eg$} if $\geog = \geodeff^2 \eg$, where $\geodeff$ is an
almost geodesic defining function of order $\reg$.

\subsubsection{Almost Geodesic Coordinates}

Given a $C^{\reg}$ asymptotically hyperbolic metric $\eg$ and an
almost geodesic compactification $\geog = \geodeff^2 \eg$ of order $\reg$
associated to $\eg$, we say a boundary adapted coordinate system
$\{\geocoord^{\alpha}\}$ is an \emph{almost geodesic coordinate system
of order $\reg$} for $\geog$ if $\geocoord^{i} \in C^{\reg+1}$ with
respect to the given smooth structure, $\geocoord^0 = \geodeff$, and
in these coordinates, $\geog_{00} = 1+ \E_\reg$ and $\geog_{i0} = \E_\reg$.
Note that the components of $\geog$ are all $C^{\reg}$ in this
coordinate system.  As for existence, we have the following:

\begin{lemma} \label{ageocoordlem}
Let $\eg$ be a $C^{\reg}$ asymptotically hyperbolic metric on a
coordinate domain $U$ containing boundary portion $D$
and let $\geog = \geodeff^2 \eg$ be an almost geodesic
compactification of order $\reg$ associated to $\eg$.  Let $\{x^i, i = 1,
\ldots, n-1\}$ be coordinates for $D$, which are $C^{\reg+1}$ compatible
with the given smooth structure.  Then there is an
almost geodesic coordinate system of order $\reg$ for $\geog$ such that
the restriction of the tangential coordinates are the given boundary
coordinates.
\end{lemma}

\begin{proof}
Proposition B.1.1 in \cite{AndChr} provides us with coordinates
$\{\geocoord^{\alpha}\}$ which are $C^{\reg+1}$ with respect
to the given smooth structure and such that $\geog_{00} = 1+ \E_\reg$ and
$\geog_{i0} = \E_\reg$.  Moreover, they can be chosen so that
$\geocoord^i = x^i$ on the boundary.

We now work in this coordinate system to show that $\geocoord^0$ and
$\geodeff$ agree modulo $\E_{\reg+1}$.  We already know that
$\geocoord^0$ agrees with $\geodeff$ on the boundary, since they are
both defining functions.  Moreover, since $\geodeff$ is an almost
geodesic defining function for $\geog$, we have
\begin{equation} \label{rhoderivs}
1-\geog^{\alpha \beta} \geodeff,_{\alpha} \geodeff,_{\beta}
= \geodeff\, \E_\reg.
\end{equation}
Evaluating this at the boundary, we find that $\geodeff,_0 = 1$.  To
determine higher derivatives at the boundary, we differentiate
\eqref{rhoderivs} $l$ times with respect to $\geocoord^0$, $1 \leq l
\leq \reg$ to get
\[
(\partial_0^l \geog^{\alpha \beta})\geodeff,_{\alpha} \geodeff,_{\beta}
   + 2 \geog^{\alpha \beta}
       (\partial_0^l \geodeff,_{\alpha}) \geodeff,_{\beta}
   + \P(\geog^{-1}, \partial_0^{l-1} \geog, \partial^{l} \geodeff)
   = (\partial_0^{l}\geodeff)\E_\reg + \cdots
         + (\partial_0 \geodeff)\E_{\reg-l+1} + \geodeff\,\E_{\reg-l},
\]
where in this case, every term in the polynomial $\P$ on the left hand
side includes at least one derivative of order $\geq 2$ of $\geodeff$.
Working inductively, and using the form of $\geog$ in the given
coordinates, we find that at the boundary, this equation reduces to
\[
\partial_0^{l+1} \geodeff = 0.
\]
Hence, $\geodeff = \geocoord^0 + \E_{\reg+1}$, and therefore the form of
the metric only changes modulo $\E_\reg$ when we look at the change of
coordinates from $\geocoord^0$ to $\geodeff$.
\end{proof}

If $\eg$ is an asymptotically hyperbolic Einstein metric and $\geog$
is an associated almost geodesic compactification, then in an almost
geodesic coordinate system for $\geog$ we get an expansion for the
tangential components of $\geog$ at the boundary similar to that for
geodesic coordinates.

\begin{prop} \label{aginageolem}
Let $2 \leq \reg \leq n-2$ and let $\eg$ be a $C^{\reg}$ asymptotically
hyperbolic Einstein metric on a coordinate domain $U \subset
\overline{M}$ with boundary portion $D \subset \partial M$.  Let
$\geog$ be an almost geodesic compactification associated to $\eg$
with boundary metric $h$,
and let $\{\geodeff, \geocoord^i\}$ be almost geodesic coordinates of
order $\reg$ for $\geog$ on $U$.  Then for $0 \leq p \leq \reg$ we have
\begin{equation} \label{aginageo}
(\partial^p \geog_{ij})|_D = \P(h^{-1}, \partial_t^p h),
\end{equation}
where the polynomials on the right hand side are the same polynomials
which appear in \eqref{geocomps}.  In particular, $(\partial_0^p
\geog_{ij})|_D = 0$ when $p$ is odd.
\end{prop}

\begin{proof}
The derivation of this expression follows the methods used in
\cite{green} to generate the similar expansion for a geodesic
compactification in geodesic coordinates.  We note that if we can
verify \eqref{aginageo} for $(\partial_0^p \geog_{ij})|_D$, then the
more general result follows immediately by taking derivatives with
respect to $\geocoord^i$ on the boundary.

Writing out the condition $\eRic = -(n-1) \eg$ under the
conformal transformation $\geog = \geodeff^2 \eg$ and focusing on the
tangential components, we find that in our almost geodesic coordinates,
\begin{eqnarray} \label{agRicexp}
\lefteqn{\geodeff\, \geog_{ij},_{00} + (2-n)\geog_{ij},_0
    - \geog^{kl} \geog_{kl},_0 \geog_{ij}} \nonumber \\
& & \mbox{} - \geodeff\, \geog^{kl} \geog_{ik},_0 \geog_{jl},_0
    + \frac{\geodeff}{2}\, \geog^{kl} \geog_{kl},_0 \geog_{ij},_0
    - 2 \geodeff\, \hat{R}ic_{ij}
    = \geodeff\, \E_{\reg-2},
\end{eqnarray}
where $\hat{R}ic$ is the Ricci tensor for the induced metric on level
sets of $\geodeff$.  This equation is the same as (2.5) in
\cite{green}, except for the $\geodeff\, \E_{\reg-2}$ term arising from
the error terms on $\geog$ and its first and second derivatives.
For $1 \leq p \leq \reg-1$, we differentiate \eqref{agRicexp} $p-1$
times with respect to $\geodeff$ to get
\begin{eqnarray} \label{agRicexpdiff}
\lefteqn{\geodeff\, \partial_0^{p+1} \geog_{ij}
         + (p-n+1)\partial_0^{p} \geog_{ij}
         - \geog^{kl}(\partial_0^p \geog_{kl})\geog_{ij}} \nonumber \\
 & &  = \geodeff\, \P(\geoinv, \partial_0^{p} \geog,
                 \partial_0^{p-2} \partial_t^2 \geog)
        + \P(\geoinv, \partial_0^{p-1} \partial_t^2\geog)
	+ \geodeff\, \E_{\reg-p-1}.
\end{eqnarray}
Setting $\geodeff$ equal to zero, tracing with respect to $h$, and
then plugging back in, we recursively solve for $(\partial_0^{p}
\geog_{ij})|_D$ to get \eqref{aginageo} for $1 \leq p \leq \reg-1$.

We cannot differentiate \eqref{agRicexp} any more than this because it
would generate derivatives of $\geog$ of order $\reg+1$.  Instead, to
derive \eqref{aginageo} in the case of $\reg$ derivatives, we set $p =
\reg-1$ in \eqref{agRicexpdiff}, divide by $\geodeff$, and evaluate the
result as a limit of a difference quotient.  To do this, we need to
know something about the derivative of order $\reg-1$.  Evaluating
\eqref{agRicexpdiff} at $\geodeff = 0$, we have
\begin{equation} \label{k-1deriv}
(\reg-n)\partial_0^{\reg-1}\geog_{ij}
   - \geog^{kl}(\partial_0^{\reg-1} \geog_{kl})\geog_{ij}
   = (\P_{\reg-1})|_D,
\end{equation}
where $\P_{\reg-1} = \P(\geoinv, \partial_0^{\reg-2} \partial_t^2\geog)$ is
the second polynomial term on the right hand side of
\eqref{agRicexpdiff}.  Now, dividing \eqref{agRicexpdiff} by $\geodeff$ and
rearranging we have
\begin{eqnarray*}
\lefteqn{\partial_0^{\reg} \geog_{ij}
         + \left(\frac{(\reg-n)\partial_0^{\reg-1} \geog_{ij}
         - \geog^{kl}(\partial_0^{\reg-1} \geog_{kl})\geog_{ij}
         - \P_{\reg-1}}{\geodeff}\right)} \hspace{2in} \\
 & = & \P(\geoinv, \partial_0^{p} \geog,
                 \partial_0^{p-2} \partial_t^2 \geog)
        + \E_{0}.
\end{eqnarray*}
We then take the limit as $\geodeff \rightarrow 0$.  The difference
quotient converges since, at the boundary, the first two terms
converge to the third term by \eqref{k-1deriv}.  Moreover, the limit
is exactly what we would have gotten if we were to differentiate each
term.  Hence, we can solve for $\partial_0^{\reg} \geog_{ij}$ as above
and when we do so, the result is the same as what we would get by
formally differentiating.

As a final note, observe that even if $\geog$ is smoother than
$C^{n-2}$, we cannot make use of \eqref{agRicexpdiff} for $p \geq n-1$
because when $p = n-1$, the second term on the left hand side is zero.
Without this term, we can no longer recursively solve for the
derivative of $\geog$.
\end{proof}

\subsubsection{Curvature of an Almost Geodesic Compactification}
\label{ageotensors}

The fact that the derivatives at the boundary of an almost geodesic
compactification $\geog$ can be expressed in terms of its boundary
metric allows us to simplify coordinate expressions at the boundary
for the Ricci tensor and scalar curvature tensor of $\geog$, along
with their covariant derivatives.

\begin{lemma} \label{geodiffRicS}
Let $2 \leq \reg \leq n-2$ and let $\geog$ be an almost geodesic
compactification of order $\reg$ associated to a $C^{\reg}$
asymptotically hyperbolic Einstein metric.  In almost geodesic
coordinates of order $\reg$ and for $0 \leq p \leq \reg-2$, we have
\[
(\geocoder^p \geoRic)|_{\partial M}
          = \P(h^{-1}, \partial_t^{p+2} h)
\]
and
\[
(\geocoder^p \geoS)|_{\partial M}
          = \P(h^{-1}, \partial_t^{p+2} h).
\]
\end{lemma}

\begin{proof}
Observe that components of $\geocoder^p \geoRic$ involve at most
$p+2$ derivatives of $\geog$.  Since we are in almost geodesic
coordinates for $\geog$, we can use \eqref{aginageo} along with the facts
that $\geog_{00} = 1 + \E_m$ and $\geog_{0i} = \E_m$ to conclude that
$(\geocoder^{p} \geoRic)|_{\partial M}$ can be expressed completely
in terms of $h$ and its first $p+2$ derivatives.  This is valid as
long we are taking few enough derivatives that we can make use of
\eqref{aginageo}.  For this, we require $p+2 \leq \reg$.  The same
analysis gives us the result for $\geoS$ and its covariant
derivatives.
\end{proof}

We will ultimately be working in coordinates that are not necessarily
almost geodesic coordinates.  The next lemma shows that we can work
with different boundary adapted coordinates without losing the
important characterization of curvature tensors and their derivatives
from the above lemma.  For this lemma, and in the future, we make use
of the inward pointing normal vector $\N$ relative to a given metric
$g$.  Note that
in boundary adapted coordinates $\{x^{\alpha}\}$, $\N$
can be written in terms of the
coordinate vectors $\partial/\partial x^{\alpha}$:
\begin{equation} \label{Normal}
\N = \frac{\grad_{g}(x^0)}{|dx^0|_{g}}
   = \frac{g^{\alpha 0}}{\sqrt{g^{00}}}
     \frac{\partial}{\partial x^{\alpha}}.
\end{equation}

\begin{lemma} \label{geodiffRicSarb}
Let $2 \leq \reg \leq n-2$ and let $\geog$ be an almost geodesic
compactification of order $\reg$ associated to a $C^{\reg}$
asymptotically hyperbolic Einstein metric.  In any arbitrary boundary
adapted coordinate system that is $C^{\reg+1}$ compatible with
the given smooth structure and for $0 \leq
p \leq \reg-2$, we have
\[
(\geocoder^p \geoRic)|_{\partial M}
          = \P(\geog^{-1}, (\geog^{00})^{-\frac{1}{2}},
               \partial_t^{p+2} h)
\]
and
\[
(\geocoder^p \geoS)|_{\partial M}
          = \P(\geog^{-1}, (\geog^{00})^{-\frac{1}{2}},
               \partial_t^{p+2} h).
\]
\end{lemma}

\begin{proof}
Let $\{x^{\alpha}\}$ be a boundary adapted coordinate system that is
$C^{\reg+1}$ compatible with the given smooth structure.  By Lemma
\ref{ageocoordlem}, there is an almost geodesic coordinate system
$\{\geocoord^{\alpha}\}$ of order $\reg$ for
$\geog$ such that $\geocoord^i|_D = x^i|_D$.  In what follows, the
metric and other tensors are being expressed in terms of
$\{x^{\alpha}\}$, so in order to use the expressions in Lemma
\ref{geodiffRicS} we need to know how to express the coordinate
vectors $\partial/\partial x^{\alpha}$ in terms of $\partial/\partial
\geocoord^{\alpha}$ at the boundary.  Note that by construction, we
already have $\partial/\partial x^{i} = \partial/\partial
\geocoord^{i}$ on $\partial M$.  Moreover, $\partial/\partial
\geocoord^{0}$ is the inward pointing normal vector $\geoN$.
By rearranging \eqref{Normal}, we can solve
for $\partial/\partial x^{0}$ in terms of $\partial/\partial
\geocoord^{\alpha}$ at the boundary:
\[
\frac{\partial}{\partial x^0}
   =  \frac{\geoN}{\sqrt{\geog^{00}}}
	- \frac{\geog^{i0}}{\geog^{00}}
	  \frac{\partial}{\partial x^i}
   =  \frac{1}{\sqrt{\geog^{00}}}
        \frac{\partial}{\partial \tilde{x}^0}
	- \frac{\geog^{i0}}{\geog^{00}}
	  \frac{\partial}{\partial \tilde{x}^i}.
\]
With this, we compute
\begin{eqnarray*}
\lefteqn{(\tilde{\nabla}^p \tilde{R}ic)
       \left(\frac{\partial}{\partial x^{\alpha}},
             \frac{\partial}{\partial x^{\beta}},
	     \frac{\partial}{\partial x^{\mu_1}},
	     \frac{\partial}{\partial x^{\mu_2}}, \ldots,
	     \frac{\partial}{\partial x^{\mu_p}}\right)} \hspace{1in} \\
& = & F^{\mu \delta \eta_1 \eta_2 \ldots\eta_p}
       _{\alpha \beta \mu_1 \mu_2 \ldots \mu_p}
           \tilde{\nabla}^p \tilde{R}ic
      \left(\frac{\partial}{\partial \tilde{x}^{\mu}},
	    \frac{\partial}{\partial \tilde{x}^{\delta}},
	    \frac{\partial}{\partial \tilde{x}^{\eta_1}},
            \frac{\partial}{\partial \tilde{x}^{\eta_2}}, \ldots,
            \frac{\partial}{\partial \tilde{x}^{\eta_p}}\right)
\end{eqnarray*}
where $F^{\eta_1 \eta_2 \cdots \eta_p}_{\mu_1 \mu_2 \cdots \mu_p} =
\P(\geoinv, \georoot)$, by the above facts about coordinate vectors.
The covariant derivative of Ricci on the right hand side can then be
replaced by the polynomial in Lemma \ref{geodiffRicS} to get the
result.
\end{proof}

Lemma \ref{geodiffRicSarb} gives us formulas for the Ricci tensor and
scalar curvature of an almost geodesic compactification at the
boundary in arbitrary boundary adapted coordinates, but they are not
explicit.  Explicit formulas for the Ricci tensor and scalar curvature
in dimension 4 are presented in
\cite{Anderson}, and the formula for scalar curvature in general
dimension is provided in \cite{Lee2}.

We finish this section with a consideration of the second fundamental
form $\geoA$ of the boundary.  Our convention here, and in the future,
is that the normal vector used to define the second fundamental form
points inward.

\begin{lemma} \label{secondfundform}
Let $\geog$ be an almost geodesic compactification for a
$C^{2}$ asymptotically hyperbolic Einstein
metric.  Then $\geoA = 0$ on $D$.
\end{lemma}

\begin{proof}
Since the second fundamental form is a tensor, we can prove this
result in coordinates of our choice, so we use almost geodesic
coordinates $\{\geocoord^{\alpha}\}$.  Since these are boundary
adapted coordinates, we have
\[
\geoA_{ij} = (\geog^{00})^{-\frac{1}{2}}\,\geoChr_{ij}^0,
\]
Writing out the Christoffel symbol and using the form of $\geog$ in
its almost geodesic coordinates, along with the fact that $(\partial_0
\geog_{ij})_D = 0$ by Proposition \ref{aginageolem}, we have the
result.
\end{proof}


\section{Deriving the Boundary Value Problem} \label{BVPchap}

In this section, we derive a boundary value
problem for components of a compactification of an asymptotically
hyperbolic Einstein metric.  The compactification we use will not be an
almost geodesic compactification, but rather one which has constant
scalar curvature.  Nonetheless, an almost geodesic compactification
will come into the picture as a tool for deriving boundary equations
for the system.

\begin{prop} \label{BVProp}
Let $\overline{M}$ be an $n$-dimensional $C^{\infty}$ manifold with
boundary, $n \geq 4$ and even.  Let $U \subset \overline{M}$ be a
coordinate domain with boundary portion $D = U \cap \partial M$.  If
$n > 4$, let $\eg$ be a $C^{n-2}$ asymptotically hyperbolic Einstein
metric on $U \cap M$, and if $n = 4$, let $\eg$ be a $C^3$ conformally
compact Einstein metric on $U \cap M$.  Let $h \in C^{n-2}(D)\
(C^3(D)\ \mathrm{if}\ n=4)$ be an element of the conformal infinity of
$\eg$ and suppose that $\eg$ has a constant scalar curvature
compactification $\csg \in C^{n-2}(U)\ (C^3(U)\ \mathrm{if}\ n=4)$
with scalar curvature $\csS = -n(n-1)$ and boundary metric $h$.
Finally, suppose that $\csg$ has boundary adapted harmonic coordinates
$\{x^{\alpha}\}$ in $C^{n-1}(U)$ ($C^4(U)$ if $n=4$).  Then, in these
harmonic coordinates, the components of $\csg$ solve a system and
accompanying boundary equations of the following form:
\begin{itemize}
\item On $U \cap M$:

\begin{equation} \label{system}
\D_{\frac{n}{2}} \csg_{\alpha \beta}
      + \P(\csg^{-1}, \partial^{n-1} \csg) = 0.
\end{equation}
\item On $D$:
\begin{itemize}
\item Equations of order 0:  For $1 \leq i,j \leq n-1$,
\begin{equation} \label{order_0}
\csg_{ij} = h_{ij}.
\end{equation}
\item Equations of order 1:  For $0 \leq \alpha \leq n-1$,
\begin{equation} \label{order_1}
\csg^{\eta \beta} \partial_{\eta} \csg_{\alpha \beta}
      - \frac{1}{2}\csg^{\eta \beta}
        \partial_{\alpha} \csg_{\eta \beta} = 0.
\end{equation}
\item Equations of order 2:  For $1 \leq i,j \leq n-1$,
\begin{equation} \label{order_2}
\D \csg_{ij} + \P(\csg^{-1}, (\csg^{00})^{-\frac{1}{2}},
                      \partial \csg, \partial_t^2 h) = 0.
\end{equation}
\item Equations of order 3:  For $0 \leq \alpha \leq n-1$,
\begin{equation} \label{order_3}
\csg^{\eta \beta} \partial_{\eta} \D \csg_{\alpha \beta}
      + \P(\csg^{-1}, \partial^2 \csg) = 0. 
\end{equation}
\item Equations of order $2l$, $2 \leq l \leq (n/2)-1$: For $0 \leq
  \alpha, \beta \leq n-1$,
\begin{equation} \label{order_2l}
\D_l \csg_{\alpha \beta}
       + \P(\csg^{-1}, (\csg^{00})^{-\frac{1}{2}},
            \partial^{2l-1}\csg, \partial_t^{2l} h) = 0.
\end{equation}
\end{itemize}
\end{itemize}

Additionally, the following formula for the derivatives of the second
fundamental form holds on $D$.  For $1 \leq k \leq n-1$,
\begin{equation} \label{difA}
\partial_k \csA_{ij} = \frac{(\csg^{00})^{\frac{1}{2}}\,\csg_{ij}}
                         {2(2-n)} \D \csg_{k0} 
                       + \P(\csinv, \csroot, \partial \csg,
                            \partial_t^2 h).
\end{equation}
\end{prop}

Before focusing on the derivation of \eqref{system}--\eqref{difA},
there are a number of comments to be made with regard to Proposition
\ref{BVProp}.  First, by Proposition \ref{intreg}, we know that $\csg$
is $C^{\infty}$ in its harmonic coordinates on the interior, so
\eqref{system} makes sense classically.  Second, the reason for the
stronger a priori regularity in dimension four is to accommodate the
boundary equations of order 3.  Third, we note that
the regularity of $\csg$ and $h$ are preserved when moving to the
given harmonic coordinates, since the coordinates are one degree
smoother than $\csg$ and $h$.  Fourth, observe that for the equations
of order greater than 3, there is no difference in behavior between
the tangential and non-tangential components of $g$.  The reason for
this will become clear in the course of the proof.  Finally, \eqref{difA} is part of the boundary value problem in that it provides important relations among the second derivatives of $\csg$.  As a consequence, it will play a crucial role in the proof of boundary regularity for the system.  In section
\ref{syschap}, we will see an alternative way to incorporate \eqref{difA} into the rest of the boundary data.

To prove Proposition \ref{BVProp}, we start by introducing the so
called Ambient Obstruction Tensor $\mathcal{O}$ and we recall that it
vanishes for any metric which is conformal to an Einstein metric.  The
equation $\mathcal{O}_{\alpha \beta} = 0$ then reduces to
\eqref{system}.  Moving on to the boundary equations, we will find
that the equations of order 1 and 3 follow as restrictions to the
boundary of equations that hold on the interior.  The remaining
equations, including \eqref{difA}, will not come as naturally.
They essentially arise from the fact that through an almost geodesic
compactification, powers of the Laplacian of the Ricci tensor can be
expressed in terms of boundary data, modulo lower order terms.

\subsection{The System via the Ambient Obstruction Tensor}

In \cite{FeffG} (Proposition 3.5), Fefferman and Graham introduced a
generalization of the Bach tensor in each even dimension
called the \emph{Ambient Obstruction Tensor} $\mathcal{O}$.  In indices, it
has the following form:
\begin{equation} \label{AOT}
\mathcal{O}_{\alpha \beta} = \frac{1}{(-2)^{\frac{n}{2}-2}(\frac{n}{2}-2)!}
                       \left(\Delta^{\frac{n}{2}-1} P_{\alpha \beta}
		       - \frac{1}{2(n-1)}
		       \Delta^{\frac{n}{2}-2} S,_{\alpha \beta}\right)
		       + \P(\csg^{-1}, \partial^{n-1}\csg)
\end{equation}
where
\[
P_{\alpha \beta} = \frac{1}{n-2}\left(Ric_{\alpha \beta}
                                  - \frac{1}{2(n-1)}S\,
				    g_{\alpha \beta}\right).
\]
When $n = 4$ this tensor is the Bach tensor, and in general, it is
symmetric, trace free, conformally invariant with weight $(n-2)/2$, and
equal to zero when the metric is Einstein.  See \cite{GHir} for details
surrounding these facts.

\begin{proof}[Proof of the System Equations]
Since we are in harmonic coordinates for a $C^2$ constant scalar curvature
compactification for an Einstein metric, $\csg \in C^{\infty}(U \cap
M)$ by Proposition \ref{intreg}, so the ambient obstruction tensor is
well defined.  Moreover, \eqref{Ric} can be used for
Ricci and $\csS$ is constant, so \eqref{AOT} reduces to
\[
\mathcal{O}_{\alpha \beta} = \frac{1}{(-2)^{\frac{n}{2}-1}(\frac{n}{2}-2)!(n-2)}
                       \left(\D_{\frac{n}{2}} \csg_{\alpha \beta}\right)
		        + \P(\csg^{-1}, \partial^{n-1} \csg).
\]
From the properties of the ambient obstruction tensor above, this is
equal to zero since $\csg$ is conformally Einstein.  After
multiplying by the leading constant, we have \eqref{system}.
\end{proof}

\subsection{Boundary Equations of Order 0, 1, and 3}

As mentioned, these equations are almost immediate.

\begin{proof}[Proof of the Boundary Equations of Order 0, 1, and 3]
For the equations of order 0, we have \eqref{order_0} since $\csg$
restricts to $h$ on $D$ and we are working in boundary adapted
coordinates.

The equations of order 1 follow from the statement that the given
coordinates are harmonic with respect to $\csg$.  The equation $\Delta
x^{\eta} = 0$ is equivalent to $\csg^{\alpha \beta} \csChr_{\alpha
\beta}^{\eta} = 0$.  Writing out the Christoffel symbol in terms of
derivatives of the metric, lowering the free index, and relabeling, we
obtain \eqref{order_1}.

Finally, the equations of order 3 follow from the contracted Bianchi
identity, the fact that $\csS$ is constant, and \eqref{Ric}.  We have
\[
0 =  \frac{1}{2}\csS,_{\alpha}
   =  \cscoder^{\beta} \csRic_{\alpha \beta}
   =  \csg^{\eta \beta}
          \left(-\frac{1}{2}(\D \csg_{\alpha \beta}),_{\eta}
	  + \P(\csg^{-1}, \partial^2 \csg)\right),
\]
and \eqref{order_3} follows.
\end{proof}

\subsection{Boundary Equations of Order $2l$ and Equation \eqref{difA}}

For these equations, we relate tensors associated with $\csg$ to
tensors associated to an almost geodesic compactification of order
$n-2$ which is conformally related to $\csg$.  The tensors associated
with $\csg$ give rise to the operator $\D_l$ acting on the components of
the metric and, through Lemma \ref{geodiffRicSarb}, the related tensors
for the almost geodesic metric can be expressed in terms of boundary
data.  The bulk of the work is in analyzing the conformal factor
relating $\csg$ and the almost geodesic compactification.

\subsubsection{Boundary Equations of Order 2 and Equation \eqref{difA}}

To generate the boundary equations of order 2, we will look at the way
that the Ricci tensor changes under a conformal transformation to an
almost geodesic metric.  This will give us a formula for $\D$ acting
on components of $\csg$, but will also introduce derivatives of
the conformal factor.  We will deal with these by analyzing the second
fundamental form and the scalar curvature under the conformal change.

By Lemma \ref{ageodefflem}, we can choose an almost geodesic
compactification $\geog$ of order $n-2$ for $\eg$ such that $\iota^*
\geog = \iota^* \csg = h$.  Then we have $\csg = \phi^2 \geog$ with
$\phi \in C^{n-2}(U)$ and $\phi|_{D} = 1$.  In this and subsequent
sections, we will be passing back and forth between these two metrics,
so we observe that when taking derivatives we have
\[
\partial^p \csg = \P(\partial^p \phi, \partial^p \geog)
\]
and
\[
\partial^p \geog = \P(\phi^{-1}, \partial^p \phi, \partial^p \csg).
\]
We start with some lemmas that tell us how to deal with derivatives of
$\phi$.

\begin{lemma} \label{phiderlem}
The first derivatives of $\phi$ at the boundary are as follows:
\[
(\phi,_i)|_D = 0
\]
and
\begin{equation} \label{phider0}
(\phi,_0)|_D = \P(\csg^{-1}, (\csg^{00})^{-1},
                  \partial \csg). 
\end{equation}
\end{lemma}

\begin{proof}
The tangential derivatives are immediate since $\phi$ is constant on
the boundary.  For the transverse derivative, we use the second
fundamental form.  Looking at the way that the
second fundamental form changes under our conformal change, we have
\[
\geoA_{ij} = \phi^{-1} \csA_{ij}
               + \phi^{-2} \phi,_{\alpha} \N^{\alpha} \csg_{ij}.
\]
We know that $\phi|_D = 1$ and $(\phi_i)|_D = 0$.  Moreover, $\geoA$
is zero by Lemma \ref{secondfundform} and $\N^{\alpha} = \csg^{\alpha
0}(\csg^{00})^{-\frac{1}{2}}$ by \eqref{Normal}.  Hence this formula
simplifies substantially to
\begin{equation} \label{simpA}
\csA_{ij} = -\phi,_0 (\csg^{00})^{\frac{1}{2}} \csg_{ij}.
\end{equation}
We are in boundary adapted coordinates so $\csA_{ij} =
(\csg^{00})^{-\frac{1}{2}}\,\csChr_{ij}^0$.  Plugging this in, along
with the fact that $\csChr_{ij}^0 = \P(\csg^{-1}, \partial \csg)$, and
solving for $\phi_0$, we get the result.
\end{proof}

\begin{lemma} \label{phiLaplem}
At the boundary,
\[
\geoLap \phi|_D = \P(\csg^{-1}, (\csg^{00})^{-\frac{1}{2}},
                  \partial \csg, \partial_t^2 h).
\]
\end{lemma}

\begin{proof}
For this formula, we consider how scalar curvature changes under the
given conformal change:
\[
S = \phi^{-2}\tilde{S} + \phi^{-3}(2-2n)\geoLap\phi
    - \phi^{-4}(n-1)(n-4)|d\phi|_{\geog}^2.
\]
Solving for $\geoLap\phi$, we have
\begin{equation} \label{Lapphi}
\geoLap\phi = \frac{1}{2(n-1)}
                     (\phi \geoS - \phi^3 \csS
                      - (n-4)(n-1)\phi^{-1}|d\phi|_{\geog}^2).
\end{equation}
Restricting to the boundary, we apply Lemma \ref{geodiffRicSarb} to
deal with $\geoS$.  This introduces first derivatives of $\geog$
which, as mentioned earlier, can be written in terms of first
derivatives of $\csg$ and first derivatives of $\phi$.  Then, we apply
Lemma \ref{phiderlem} to deal with the first derivatives of $\phi$,
and we replace $\geog$ with $\csg$ since they are equal on $D$.
Finally, $\csS$ is constant, so the result follows.
\end{proof}

To derive the boundary equations of order 2 and \eqref{difA}, we will
look at how various components of the Ricci tensor change under the
given conformal change.  For general components, we have
\begin{eqnarray} \label{conRic}
\csRic_{\alpha \beta} & = & \tilde{R}ic_{\alpha \beta}
  + \phi^{-1}\Bigl((2-n)\geocoder_{\alpha}\geocoder_{\beta}\phi
              - \geoLap\phi\,\geog_{\alpha \beta}\Bigr) \nonumber \\
 &   & \mbox{} + \phi^{-2}\Bigl((3-n)|d\phi|^2_{\geog}\,\geog_{\alpha \beta}
               + 2(n-2)\phi,_{\alpha} \phi,_{\beta}\Bigr).
\end{eqnarray}
With this and the lemmas above, we are ready to derive the boundary
equations of order 2 and \eqref{difA}.

\begin{proof}[Proof of the Boundary Equations of Order 2]
When restricted to $D$ the tangential components of
\eqref{conRic} reduce to
\begin{equation} \label{conRicij}
\csRic_{ij} = \tilde{R}ic_{ij}
  + (2-n)\geocoder_{i}\geocoder_{j}\phi
  - \geoLap\phi\geog_{ij}
  + (3-n) |d\phi|^2_{\geog}\geog_{ij}
  + 2(n-2) \phi,_{i} \phi,_{j}.
\end{equation}
Now, making a number of substitutions, we will arrive at
\eqref{order_2}.  For the left hand side, we again use \eqref{Ric}.
For the first term on the right hand side, we use Lemma
\ref{geodiffRicSarb}.  Then, as in Lemma \ref{phiLaplem}, we replace
the first derivatives of $\geog$ with first derivatives of $\csg$ and
first derivatives of $\phi$, and then we use Lemma \ref{phiderlem} to
handle the first derivatives of $\phi$.

For the second term, since $\geoA$ is zero, we also have
$\geoChr_{ij}^0 = (\geog^{00})^{\frac{1}{2}} \geoA_{ij} = 0$.  Hence,
since the tangential
derivatives of $\phi$ are zero on $D$, we have
\[
\geocoder_i \geocoder_j \phi = \phi,_{ij}
                               - \phi,_{\alpha} \geoChr^{\alpha}_{ij}
                             = 0.
\]

For the third term, we use the formula from Lemma \ref{phiLaplem}, and
then substitute as above to eliminate the first derivatives of
$\geog$.

The fourth term, which involves only $\geog$ and first derivatives of
$\phi$, is taken care of by making the same substitutions as were made
earlier. Finally, the fifth term vanishes.  Putting all this into
\eqref{conRicij}, and rearranging, we have \eqref{order_2}.
\end{proof}

\begin{proof}[Proof of Equation \eqref{difA}]
Differentiating \eqref{simpA} and using \eqref{phider0}, we have
\begin{equation} \label{Awrtphi}
\partial_k \csA_{ij} = - \phi,_{0k}(\csg^{00})^{\frac{1}{2}} \csg_{ij}
                       + \P(\csinv, \csroot, \partial \csg).
\end{equation}
To gain information about $\phi,_{0k}$, we set $\alpha = k$ and $\beta
= 0$ in \eqref{conRic} and solve for $\geocoder_k \geocoder_0 \phi$
while using various results on the other terms in \eqref{conRic}.
In particular, we use \eqref{Ric} for the left hand side, Lemma
\ref{geodiffRicSarb} for $\geoRic$, Lemma \ref{phiLaplem} for the
Laplacian of $\phi$, and Lemma \ref{phiderlem} for the first derivatives
of $\phi$ to end up with
\[
\geocoder_k \geocoder_0 \phi = \frac{-1}{2(2-n)}\D \csg_{k0}
                               + \P(\csinv, \csroot, \partial \csg,
                                    \partial_t^2 h).
\]
On the other hand, $\geocoder_k \geocoder_0 \phi = \phi,_{0k} -
\phi,_{\alpha}\geoChr_{0k}^{\alpha}$ and we can handle the first
derivatives of $\phi$ and the Christoffel symbol using Lemma
\ref{phiderlem} and the fact that $\partial \geog = \P(\phi^{-1},
\partial \phi, \partial \csg)$.  Substituting this and rearranging, we
get
\[
\phi,_{0k} = \frac{-1}{2(2-n)}\D \csg_{k0}
                               + \P(\csinv, \csroot, \partial \csg,
                                    \partial_t^2 h).
\]
Finally, substituting this into \eqref{Awrtphi} we get \eqref{difA}.
\end{proof}

\subsubsection{Boundary Equations of Order $2l$, $2 \leq l \leq (n/2)-1$}

Generating the boundary equations for higher even orders follows the
same basic idea as for the boundary equations of order 2.  Here,
looking at how the $l$th power of the Laplacian of Ricci changes under
a conformal change of the metric gives us an expression for $\D_l$
applied to components of the metric.  However it also introduces
derivatives of the conformal factor, so to start we will determine how
to keep track of derivatives of the conformal factor inductively in
order to get the result.  In the process, we will see why we do not
need to distinguish between tangential and non-tangential components
for these equations.

As before, let $\geog$ be an almost geodesic compactification of order
$n-2$ for $\eg$ with
$\iota^*\geog = \iota^*\csg = h$ so that $\csg = \phi^2 \geog$ with
$\phi|_D = 1$.  We begin with a formula for derivatives of $\phi$ that
we will use a number of times.

\begin{lemma}
On $U$,
\begin{equation} \label{nablapdelta}
\geocoder^p \geoLap\phi
                     =  \P(\geoinv,
                           \partial^{p}\geog,
			   \phi^{-1},
		           \partial^{p+1}\phi,
                           \geocoder^{p}\geoS).
\end{equation}
\end{lemma}

\begin{proof}
Taking covariant derivatives of \eqref{Lapphi} with respect to $\geog$
and using the fact that $\csS$ is constant, we have the result.
\end{proof}

The first application of this lemma is the following:

\begin{lemma} \label{phidercalclemma}
For $0 \leq p \leq n-4$,
\begin{equation} \label{phidercalc}
(\partial_0^{p+2} \phi)|_D
   = -\frac{1}{\csg^{00}}
     (2 \csg^{i0}\partial_0^{p+1}\partial_i\phi
          + \csg^{ij}\partial_0^{p}\partial_i\partial_j\phi)
     + \P(\csinv, \csroot, \partial^{p+1}\csg,
          \partial^{p+1}\phi, \partial_t^{p+2}h).
\end{equation}
\end{lemma}

\begin{proof}
Writing $\geocoder^p\geoLap\phi$
out directly, we have
\[
\geocoder_{\mu_1} \cdots \geocoder_{\mu_p} \geoLap\phi
   = \geog^{00}\phi,_{00\mu_1\cdots\mu_p}
     + 2 \geog^{i0}\phi,_{i0\mu_1\cdots\mu_p}
     + \geog^{ij}\phi,_{ij\mu_1\cdots\mu_p}
     + \P(\geoinv, \partial^{p+1}\phi,\partial^{p+1}\geog).
\]
Applying \eqref{nablapdelta} to the left hand side and then
rearranging, we get
\[
\phi,_{00\mu_1\cdots\mu_p}
   = -\frac{1}{\geog^{00}}
       (2 \geog^{i0}\phi,_{i0\mu_1\cdots\mu_p}
           + \geog^{ij}\phi,_{ij\mu_1\cdots\mu_p})
     + \P(\geoinv, (\geog^{00})^{-1}, \partial^{p+1}\geog,
          \phi^{-1}, \partial^{p+1}\phi, \geocoder^p \geoS).
\]
From here, we express objects associated with $\geog$ in terms of
objects associated with $\csg$ by using the fact that $\partial^m
\geog = \P(\phi^{-1}, \partial^m \phi, \partial^m \csg)$.  Also, for
$p \leq n-4$, we may use Lemma \ref{geodiffRicSarb} for the covariant
derivatives of $\geoS$.  Making these substitutions we conclude
\[
\phi,_{00\mu_1\cdots\mu_p}|_D
   = -\frac{1}{\csg^{00}}
     (2 \csg^{i0}\phi,_{i0\mu_1\cdots\mu_p}
          + \csg^{ij}\phi,_{ij\mu_1\cdots\mu_p})
     + \P(\csinv, \csroot, \partial^{p+1}\csg,
          \partial^{p+1}\phi, \partial_t^{p+2}h).
\]
In particular, this is true when all $\mu_i = 0$.
\end{proof}

We can use this result together with an inductive argument to say
even more.

\begin{lemma} \label{derphilemma}
For $0 \leq q \leq n-2$ and $0 \leq s \leq n-2-q$,
\begin{equation} \label{derphi}
(\partial_t^s \partial_0^{q}\phi)|_D
            = \P(\csinv, \csroot,
                 \partial^{q+s} \csg, \partial_t^{q+s} h).
\end{equation}
\end{lemma}

\begin{proof}
If it is already known for some $q$ that $(\partial_0^{q}\phi)|_D =
\P(\csinv, \csroot, \partial^{q} \csg, \partial_t^{q} h)$, then
differentiating tangentially proves \eqref{derphi} for $0 \leq s \leq
n-2-q$.  Hence it is enough to show the result when $s = 0$.  For
this, we use induction on $q$.  First, we have some base cases.  For
$q = 0$, we know $\phi|_{\partial M} = 1$, while for $q = 1$,
\eqref{phider0} says $\phi,_0|_D = \P(\csinv, (\csg^{00})^{-1},
\partial \csg)$.

Now let $ q \geq 2$, suppose this lemma is true for all $m < q$.
Letting $q = p+2$, we find that every derivative of $\phi$ on the
right hand side of \eqref{phidercalc} is taken care of by the
induction hypothesis.
\end{proof}

If we know that two of the derivatives on $\phi$ are actually coming
from the Laplacian, we can do a little bit better.

\begin{lemma} \label{coderplapphilemma}
For $0 \leq p \leq n-4$,
\begin{equation} \label{coderplapphi}
(\geocoder^p \geoLap \phi)|_D
            = \P(\csinv, \csroot, \partial^{p+1} \csg,
                 \partial_t^{p+2} h).
\end{equation}
\end{lemma}

\begin{proof}
By \eqref{nablapdelta}, $\geocoder^p\geoLap\phi = \P(\geoinv,
\partial^{p}\geog, \phi^{-1}, \partial^{p+1}\phi,
\geocoder^{p}\geoS)$, and as before, we use the fact that
$\partial^{m}\geog = \P(\phi^{-1}, \partial^{m} \phi,
\partial^{m}\csg)$ to eliminate derivatives of $\geog$.  Then Lemma
\ref{derphilemma} tells us how to deal with the derivatives of $\phi$,
and Lemma \ref{geodiffRicSarb} tells us how to deal with the
derivatives of $\tilde{S}$.  Inserting these equations into
\eqref{nablapdelta} produces the result.
\end{proof}

The next lemma is important for showing that we will not need
generalizations of the Neumann data for our system.

\begin{lemma} \label{derlapderphilemma}
For $0 \leq p \leq n-4$,
\begin{equation} \label{derlapderphi}
(\geocoder^{p-2} \geoLap \geocoder^{2}\phi)|_D
            = \P(\csinv, \csroot, \partial^{p+1} \csg,
                 \partial_t^{p+2} h).
\end{equation}
\end{lemma}

\begin{proof}
First note that by rearranging covariant derivatives, we have
\[
\geoLap \geocoder^{2}\phi
         = \geocoder^2 \geoLap \phi
           + \P(\geoinv, \partial^3 \geog, \partial^3 \phi).
\]
Hence,
\[
\geocoder^{p-2}\geoLap \geocoder^{2}\phi
          = \geocoder^p\geoLap\phi
	    + \P(\geoinv, \partial^{p+1} \geog,
                 \partial^{p+1} \phi).
\]
Once again, we have $\partial^{m}\geog = \P(\phi^{-1}, \partial^{m}
\phi, \partial^{m}\csg)$.  When we restrict to the boundary, Lemma
\ref{derphilemma} and Lemma \ref{coderplapphilemma} allow us to
replace all derivatives of $\phi$ and this gives us the result.
\end{proof}

With the derivatives of $\phi$ understood, we are ready to finish the
proof of Proposition \ref{BVProp} by deriving the boundary equations
of order $2l$.

\begin{proof}
[Proof of the Boundary Equations of Order $2l$, $2 \leq l \leq (n/2)-1$]
To start, we look at how the covariant derivatives of the Ricci tensor
are affected by a conformal change.  For $m = l-1$,
the conformal change formula for covariant derivatives gives us
\[
\nabla_{\mu_1}\cdots\nabla_{\mu_{2m}}Ric_{\alpha \beta}
    = \geocoder_{\mu_1}\cdots\geocoder_{\mu_{2m}}
      Ric_{\alpha \beta}
      + \P(\csinv, \partial^{2m}\csg, \phi^{-1},\partial^{2m}\phi,
           \partial^{2m-1}Ric).
\]
Next, we use the conformal change formula for the Ricci tensor to get
\begin{eqnarray*}
\nabla_{\mu_1}\cdots\nabla_{\mu_{2m}}\csRic_{\alpha \beta}
  & = & \geocoder_{\mu_1}\cdots\geocoder_{\mu_{2m}}
        \tilde{R}ic_{\alpha \beta} \\
  &   & + \phi^{-1}\Bigl((2-n)\geocoder_{\mu_1}\cdots\geocoder_{\mu_{2m}}
                       \geocoder_{\alpha}\geocoder_{\beta}\phi
                       - \geocoder_{\mu_1}\cdots\geocoder_{\mu_{2m}}
                         \geoLap\phi\,\geog_{\alpha \beta}\Bigr) \\
  &   & + \P(\csinv, \partial^{2m+1}\csg, \phi^{-1}, \partial^{2m+1}\phi).
\end{eqnarray*}
Note that the polynomial representing lower order terms loses its
dependence on $Ric$, with the slack being picked up by the extra
derivative on $g$.  We also pick up an extra derivative of $\phi$ from
the conformal change of the Ricci tensor.

Now, we trace with respect to $g = \phi^2 \geog$ to get
\begin{eqnarray*}
\Delta^m \csRic_{\alpha \beta}
   & = & \phi^{-2m}\geoLap^m\tilde{R}ic_{\alpha\beta} \\
   &   & + \phi^{-2m-1}\Bigl((2-n)\geoLap^m
                        \geocoder_{\alpha}\geocoder_{\beta}\phi
		   - \geoLap^m\geoLap\phi\,
                     \geog_{\alpha\beta}\Bigr) \\
   &   & + \P(\csinv, \partial^{2m+1}\csg, \phi^{-1}, \partial^{2m+1}\phi).
\end{eqnarray*}
Restricting to the boundary, $\phi = 1$ and we may apply all our
lemmas to various terms on the right hand side.  In particular, we use
Lemma \ref{geodiffRicSarb} for the first term, Lemma
\ref{derlapderphilemma} for the second term, Lemma
\ref{coderplapphilemma} for the third term, and Lemma
\ref{derphilemma} for the terms in the polynomial to conclude
\[
\Delta^m \csRic_{\alpha \beta}
    = \P(\csinv, \csroot, \partial^{2m+1}\csg,\partial_t^{2m+2}h).
\]
Note that we cannot draw this conclusion when $m = 0$, because then
there are no Laplacians to rearrange on the Hessian of $\phi$ in the
second term, and so we cannot make use of Lemma
\ref{derlapderphilemma}.  This is why we need the boundary equations
of order 3.

The last step is to use \eqref{Ric} on the left hand side to get
\eqref{order_2l}.
\end{proof}


\section{Local and Global Regularity} \label{regularitychap}

Using the boundary value problem that was derived in the last section,
we are ready to prove Theorems \ref{mainlocalthm} and
\ref{mainglobalthm}, which we restate here.

\medskip
\noindent \textbf{Theorem \ref{mainlocalthm}.}
\emph{Let $\overline{M}$ be an $n$-dimensional $C^{\infty}$ manifold with
boundary, $n \geq 4$ and even.  Let $p \in \partial M$ and let $U
\subset \overline{M}$ be a neighborhood of $p$ with boundary portion
$D = U \cap \partial M$.  For $\basereg \geq n$ and $0 < \holdreg < 1$, let $\eg$ be a
$C^{\basereg-1,\holdreg}$ conformally compact Einstein metric on $U \cap M$.
Suppose that the conformal infinity of $\eg$ contains a metric $h \in
C^{\newreg,\newholdreg}(D)$, where $\newreg \geq \basereg$ and $0 <
\newholdreg < 1$.  Given a $C^{\infty}$ coordinate system on the
boundary in a neighborhood of $p$, there is a coordinate system on a
neighborhood $V \subset U$ of $p$ that is $C^{\basereg, \holdreg}$ compatible
with the given smooth structure and which restricts to the given
coordinates on the boundary, and there is a defining function $\rho
\in C^{\basereg-1,\holdreg}(V)$ in the new coordinates, such that $\rho^2
\eg$ has boundary metric $h$ and in the new coordinates, $\rho^2 \eg
\in C^{\newreg,\newholdreg}(V)$.
}

\medskip
\noindent \textbf{Theorem \ref{mainglobalthm}.}
\emph{Let $\overline{M}$ be a compact $n$-dimensional $C^{\infty}$ manifold
with boundary, $n \geq 4$ and even.  For $\basereg \geq n$ and $0 < \holdreg < 1$,
let $\eg$ be a $C^{\basereg-1,\holdreg}$ conformally compact Einstein metric on
$M$.  Suppose that the conformal infinity of $\eg$ contains a metric
$h \in C^{\newreg,\newholdreg}(\partial M)$, where $\newreg \geq \basereg$
and $0 < \newholdreg < 1$.  Then there is a
 $C^{\basereg, \holdreg}$ diffeomorphism
$\Psi: \overline{M} \longrightarrow \overline{M}$ which restricts to the identity on the boundary and there is a defining function $\check{\rho} \in C^{\basereg-1,\holdreg}(\overline{M})$ such that
$\check{\rho}^2\Psi^*(\eg)$ has boundary metric $h$ and
$\check{\rho}^2 \Psi^*(\eg) \in C^{\newreg,\newholdreg}(\overline{M})$.
}

\medskip
To prove Theorem \ref{mainlocalthm}, we apply classical elliptic
regularity results to the boundary value problem that was derived in
the previous section.  We use a bootstrap argument, each step of which
requires its own inductive argument.  Theorem \ref{mainglobalthm}
follows from Theorem \ref{mainlocalthm} by studying the regularity of
the atlas of harmonic coordinates that we use and then applying an approximation theorem.

\subsection{Local Boundary Regularity}

Before proving Theorem \ref{mainlocalthm}, we have two supporting
lemmas.  The first lemma shows that the various types of boundary
conditions all produce the same regularity result, and is basically an
application of classical second order elliptic boundary regularity
results.  In order to apply these results, we observe that when
composing $\D$ with $\D_{l}$, we have
\begin{equation*}
\D(\D_l g_{\alpha \beta}) = \D_{l+1} g_{\alpha \beta}
                            + \P(g^{-1}, \partial^{2l+1} g),
\end{equation*}
where the polynomial of lower order terms is zero if $l = 0$.
Hence, we can think of $\D_{l+1} g_{\alpha \beta}$ as a second order
operator acting on $\D_l g_{\alpha \beta}$, modulo lower order terms.

\begin{lemma} \label{bvalreglem}
Let $U \subset \overline{M}$ be a boundary adapted coordinate domain
with boundary portion $D$, let $\csg$ be a metric on $U$, and let $h$
be a metric on $D$.  Let $p \geq 0$, $0 \leq l \leq (n/2) - 1$,
and $0 < \holdreg < 1$.  Suppose
$\csg_{\alpha \beta} \in C^{p+2l+1,\holdreg}(U) \cap C^{2l+2}(U \cap
M)$ and $h_{ij} \in C^{p+2l+2,\holdreg}(D)$.  Also, suppose
$\D_{l+1}\csg_{\alpha \beta} \in C^{p,\holdreg}(U)$.  Moreover,
referring to \eqref{order_0}--\eqref{order_2l} and \eqref{difA} of Proposition
\ref{BVProp}, suppose the components of $\csg$ and $h$ solve the
boundary equations of order 0 and 1 and \eqref{difA} if $l = 0$, order
2 and 3 if $l = 1$, and order $2l$ if $l \geq 2$.  Then $\D_l
\csg_{\alpha \beta} \in C^{p+2,\holdreg}(U)$.
\end{lemma}

\begin{proof}
Note that we may use regularity results for linear equations
since we are working with a fixed metric.  With this in mind, we
consider three cases.

\textit{Case 1: $2 \leq l \leq (n/2)-1$.}  By our observation about
composition of $\D$ with $\D_l$, we have
\begin{equation} \label{linteq}
\D(\D_l \csg_{\alpha \beta})
    = \D_{l+1}\csg_{\alpha \beta} + \P(\csinv, \partial^{2l+1} \csg)
\end{equation}
on $U \cap M$.  Also, as given by \eqref{order_2l}, the boundary
equation of order $2l$ is
\begin{equation} \label{lbdyeq}
(\D_l \csg_{\alpha \beta})|_D
   = \P(\csinv, \csroot, \partial^{2l-1} \csg, \partial_t^{2l} h)
\end{equation}
We view \eqref{linteq} as a linear scalar equation for $\D_l
\csg_{\alpha \beta}$ and \eqref{lbdyeq} as a Dirichlet boundary
condition.  By the regularity hypotheses provided, we have $\D_l
\csg_{\alpha \beta} \in C^{p+1,\holdreg}(U) \cap C^{2}(U \cap M)
\subset C^0(U) \cap C^2(U \cap M)$.
Moreover, the right hand side of \eqref{linteq} is in
$C^{p,\holdreg}(U)$, and the right hand side of \eqref{lbdyeq} is in
$C^{p+2,\holdreg}(D)$.  Hence, by local regularity results for the
Dirichlet problem \cite{G&T}, we may conclude that $\D_l \csg_{\alpha
\beta} \in C^{p+2,\holdreg}(U)$.

\textit{Case 2: $l = 1$.}  In this case, similar to \eqref{linteq}, we
have
\begin{equation} \label{oneinteq}
\D(\D \csg_{\alpha \beta})
      = \D_2\csg_{\alpha \beta} + \P(\csinv, \partial^3 \csg)
\end{equation}
on $U \cap M$, but now our boundary conditions require
more care than the previous argument.  We proceed in three steps.  In
the first step, we focus only on tangential components of $\csg$.  For
these, the boundary equations of order 2 provide Dirichlet conditions
and combined with \eqref{oneinteq}, we conclude that
$\D\csg_{ij} \in C^{p+2,\holdreg}(U)$ just as in case 1.

For the next two steps, we will use the boundary equations of order
3 by isolating the terms with $\beta = 0$ in \eqref{order_3}.
Writing this equation out and
dividing everything by $(\csg^{00})^{\frac{1}{2}}$, we have
\begin{equation} \label{normalu}
\N(\D \csg_{\alpha 0})
    = - \csroot \csg^{\eta i} \partial_{\eta} (\D \csg_{\alpha i})
      + \P(\csinv, \csroot, \partial^2 \csg).
\end{equation}

For the second step, we let $\alpha = j$.  Then the first term on the
right hand side of \eqref{normalu} is in $C^{p+1,\holdreg}(D)$ as a
result of the first step above.  The second term on the right is also
in $C^{p+1,\holdreg}(D)$.  Moreover, $\D \csg_{\alpha
\beta} \in C^{1}(U) \cap C^{2}(U \cap M)$ and the right hand side of
\eqref{oneinteq} is in $C^{p,\holdreg}(U)$.  Hence, treating the pair
\eqref{oneinteq}, \eqref{normalu} as a
Neumann problem for a linear scalar equation, we may use elliptic
regularity \cite{Miranda} to conclude $\D \csg_{0j} \in
C^{p+2,\holdreg}(U)$.

For the third step we let $\alpha = 0$ in \eqref{normalu}.  Then, by
the second step, the first term on the right hand side is in
$C^{p+1,\holdreg}(D)$, as is the second term.  Therefore, by the
same argument as in the second step, we have $\D \csg_{00} \in
C^{p+2,\holdreg}(U)$.  These three steps together handle every
component of $\D \csg$.

\textit{Case 3: $l = 0$.}  This case is similar to case 2 in that we
deal with the tangential and nontangential components of $\csg$
separately.  For any components we have
\begin{equation} \label{0inteq}
\D \csg_{\alpha \beta}
\in C^{p,\holdreg}(U),
\end{equation}
which is provided for us by the hypotheses.  As in case 2, we work in
three steps, although the second and third steps are reversed relative
to case 2.  For the first step, we focus on the tangential components
and the argument is essentially the same as for the tangential
components in case 2.  We are given that $\csg_{ij} \in C^{0}(U) \cap
C^2(U \cap M)$ and the boundary equations of order 0 are just the
statements that $\csg_{ij} = h_{ij}$, which are given to be in
$C^{p+2,\holdreg}(D)$.  Hence, by boundary regularity results for the
Dirichlet problem, we have $\csg_{ij} \in C^{p+2,\holdreg}(U)$.

For the second and third steps, we will make use of the boundary
equations of order 1, but first an analysis of \eqref{difA} will
provide an essential simplification.  To start, we note that $\D
\csg_{k0}$, the metric $\csg$ and all first derivatives of $\csg$, and up to
second derivatives of $h$ are in $C^{p,\holdreg}(U)$.  Hence
\eqref{difA} tells us that $\partial_k \csA_{ij} \in
C^{p,\holdreg}(D)$ for $1 \leq k \leq n-1$ and so $\csA_{ij} \in
C^{p+1,\holdreg}(D)$.  Since we are in boundary adapted coordinates,
\begin{eqnarray} \label{Abdy}
\csA_{ij} & = & \csroot \csChr_{ij}^0 \nonumber \\
          & = & \frac{1}{2}\csroot \csg^{0 \alpha}
                (\partial_j \csg_{i \alpha}
                 + \partial_i \csg_{j \alpha}
                 - \partial_{\alpha} \csg_{ij}).
\end{eqnarray}
From step 1, all tangential components of $\csg$ are in
$C^{p+2,\holdreg}(U)$.  With this and the regularity determined for
$\csA_{ij}$, we rearrange \eqref{Abdy} and conclude
\begin{equation} \label{symgder}
\csg_{i0},_j + \csg_{j0},_i \in C^{p+1,\holdreg}(D).
\end{equation}
Multiplying by $\csg^{ij}$ and summing, we also have
\begin{equation} \label{trgder}
\csg^{ij} \csg_{i0},_j \in C^{p+1,\holdreg}(D).
\end{equation}

With these facts, we focus on the boundary equations of order 1.
For the second step, we look at \eqref{order_1} with $\alpha = 0$.
Writing this out by separating terms with 0 as an index, we have
\[
\csg^{ij} \partial_i \csg_{0j} + \csg^{0j} \partial_0 \csg_{0j}
  + \csg^{j0} \partial_j \csg_{00} + \csg^{00} \partial_0 \csg_{00}
  - \frac{1}{2}(\csg^{ij} \partial_0 \csg_{ij}
                + 2 \csg^{0j} \partial_0 \csg_{0j}
                + \csg^{00} \partial_0 \csg_{00}) = 0.
\]
The second term cancels with the middle term in parentheses and the
fourth term partially cancels the last term in parentheses.
Rearranging then gives
\[
(\csg^{j0} \partial_j + \frac{1}{2}\csg^{00}\partial_0)\csg_{00}
   = \frac{1}{2}\csg^{ij} \partial_0 \csg_{ij}
     - \csg^{ij} \partial_i \csg_{0j}.
\]
By \eqref{trgder} and the fact that $\csg_{ij} \in
   C^{p+2,\holdreg}(U)$, we have
\begin{equation} \label{bdyg00}
(\csg^{j0} \partial_j + \frac{1}{2}\csg^{00}\partial_0)\csg_{00}
     \in C^{p+1,\holdreg}(D).
\end{equation}
Combined with \eqref{0inteq}, this gives us a regular oblique
derivative problem.  Since $\csg \in C^1(U) \cap C^{2}(U \cap
M)$, we apply boundary regularity for such a problem \cite{Miranda} to
conclude $\csg_{00} \in C^{p+2,\holdreg}(U)$.

For the third step, we repeat this analysis using \eqref{order_1} with
$\alpha = j$.  In this case we have
\begin{equation} \label{order_1j}
\csg^{\eta k} \partial_{\eta} \csg_{jk}
  + \csg^{i0} \partial_i \csg_{j0} + \csg^{00} \partial_0 \csg_{j0}
  - \frac{1}{2}(\csg^{ik} \partial_j \csg_{ik}
                + 2 \csg^{i0} \partial_j \csg_{i0}
                + \csg^{00} \partial_j \csg_{00}) = 0.
\end{equation}
The first term along with the first and last terms in parentheses are
in $C^{p+1,\holdreg}(U)$ by the first and second steps.  Moreover, by
\eqref{symgder} we may replace $\csg^{i0}\partial_j \csg_{i0}$ by
$-\csg^{i0}\partial_i \csg_{j0}$ modulo a term in $C^{p+1,\holdreg}$.
Using these facts, \eqref{order_1j} simplifies to
\begin{equation} \label{bdygj0}
(\csg^{i0} \partial_i + \frac{1}{2}\csg^{00}\partial_0)\csg_{j0}
     \in C^{p+1,\holdreg}(D).
\end{equation}
Therefore, by the same argument as in the second step, $\csg_{j0} \in
C^{p+2,\holdreg}(U)$.  These three steps together handle every
component of $\csg$, and the lemma is proved.
\end{proof}

Our second lemma uses Lemma \ref{bvalreglem} to provide us with each
step in the eventual bootstrap.

\begin{lemma} \label{bootstrapstep}
Let $U \subset \overline{M}$ be a boundary adapted coordinate domain
with boundary portion $D$.  Let $\reg \geq n-1$, let $\csg$ be a
metric in $C^{\reg,\holdreg}(U) \cap C^n(U \cap M)$, and let $h$ be a
metric in $C^{\reg+1,\holdreg}(D)$.  Suppose the components of $\csg$
and $h$ solve the boundary value problem
\eqref{system}--\eqref{order_2l} and \eqref{difA}.  Then the
components of $\csg$ are in $C^{\reg+1,\holdreg}(U)$.
\end{lemma}

\begin{proof}
We prove this using a reverse induction
argument on the power of $\D$ acting on components of $\csg$.  For the
base case, observe that the second term on the left hand side of \eqref{system} is in
$C^{\reg-n+1,\holdreg}(U)$, so setting $l = (n/2)-1$ and $p = \reg-n+1$, we
find that the conditions of Lemma \ref{bvalreglem} are satisfied.
Hence we may conclude that $\D_{\frac{n}{2}-1} \csg_{\alpha \beta} \in
C^{\reg-n+3,\holdreg}(U)$.

For the induction step, suppose that for some $l$ with $ 0 \leq l \leq
(n/2)-1$, we have $\D_{l+1}\csg_{\alpha \beta} \in
C^{\reg-2l-1,\holdreg}(U)$.  Setting $p = \reg-2l-1$, we find that the
conditions of Lemma \ref{bvalreglem} are satisfied, and so $\D_l
\csg_{\alpha \beta} \in C^{\reg-2l+1,\holdreg}(U)$.

The induction terminates at $l = 0$, in which case $p = \reg-1$
and we are left with $\csg_{\alpha \beta} \in C^{\reg+1,\holdreg}(U)$,
which was our goal.
\end{proof}

\begin{proof}[Proof of Theorem \ref{mainlocalthm}]
By Proposition \ref{csprop}, $\eg$ has a constant scalar curvature
compactification $\csg = \csdeff^2 \eg$ in $C^{\basereg-1,\holdreg}$ near $p$
with boundary metric $h$, where for any smooth defining function
$\hat{\rho}$, we have $\csdeff/\hat{\rho} \in C^{\basereg-1,\holdreg}$ near
$p$.  By Proposition \ref{hcoordsprop} there are harmonic coordinates
for $\csg$ in a neighborhood $V$ about $p$ which are $C^{\basereg,\holdreg}$
compatible with the given smooth structure on $U$ and which restrict
to coordinates in the given smooth structure on $D$.  Note that in
these new coordinates, $\csdeff \in C^{\basereg-1,\holdreg}(V)$.  With this
compactification and these coordinates, the conditions of Proposition
\ref{BVProp} are satisfied, so $\csg$ satisfies the resulting boundary
value problem.  We will apply Lemma \ref{bootstrapstep} inductively to
conclude the desired regularity.

We consider two cases depending on the relationship between
$\newholdreg$ and $\holdreg$.  If $\newholdreg \leq \holdreg$, we may
replace $\holdreg$ by $\newholdreg$.  Then, applying Lemma
\ref{bootstrapstep} inductively, we have the result.

On the other hand, if $\newholdreg > \holdreg$, then after applying Lemma
\ref{bootstrapstep} once, we have $\csg \in C^{\basereg,\holdreg} \subset
C^{\basereg-1,\newholdreg}$.  We then apply Lemma \ref{bootstrapstep} inductively,
starting with $\csg \in C^{\basereg-1,\newholdreg}$, and the result again follows.
\end{proof}

\subsection{Global Regularity}

Theorem \ref{mainglobalthm} now follows from Theorem
\ref{mainlocalthm}, a patching argument, and an approximation theorem.  Twice in the patching argument,
we will use the general fact, as discussed in Section \ref{harmonic},
that if a metric $g$ is $C^{k,\gamma}$ with respect to a given
coordinate chart, then harmonic coordinates are $C^{k+1,\gamma}$ with
respect to the given coordinates.

\begin{proof}[Proof of Theorem \ref{mainglobalthm}]
First observe that by the remark after Proposition \ref{csprop}, we
have a global constant scalar curvature compactification $\csg$.  Letting $\A$ be the maximal $C^{\infty}$ atlas for $\overline{M}$,  we
construct a new atlas $\B$, about each point in $\overline{M}$ by choosing
harmonic coordinates for $\csg$.  In particular, for any point in
$\partial M$, we use the coordinates in the proof of Theorem
\ref{mainlocalthm}.  Since $\csg$ is in $C^{\basereg-1,\holdreg}$ with
respect to $\A$, such coordinates are
$C^{\basereg,\holdreg}$ compatible with $\A$ since they
are harmonic.  Moreover, by Theorem \ref{mainlocalthm}, $\csg$ is in
$C^{\newreg,\newholdreg}$ in each new coordinate chart, while on the
interior, by Proposition \ref{intreg}, $\csg$ is $C^{\infty}$.  Hence,
$\csg$ is in $C^{\newreg,\newholdreg}$ on all such charts.  As a
consequence, the collection of these harmonic coordinate charts must
be $C^{\newreg+1,\newholdreg}$ compatible with one another.

Now let $id_1: (\overline{M}, \A) \longrightarrow (\overline{M}, \B)$ be the identity map.  Let
$id_2 = id_1^{-1}$ and note that $id_2^*(\eg)$ is just $\eg$ in the new atlas $\B$.  By Whitney approximation, $id_1$, which is $C^{\basereg,\holdreg}$,
can be approximated in $C^1$ by a diffeomorphism
$\Theta:(\overline{M}, \A) \longrightarrow (\overline{M}, \B)$ which is $C^{\newreg+1,\newholdreg}$ and restricts to the identity on $\partial M$.  With this, let $\Psi = id_2 \circ \Theta$ and
$\check{\rho} = \Theta^*\rho$.  By Theorem \ref{mainlocalthm}, $\rho^2 id_2^*(\eg)$ has the desired boundary characteristics.  Since $\rho^2 id_2^*(\eg)$ is in $C^{\newreg,\newholdreg}(\overline{M},\B)$, 
we have
\[
\check{\rho}^2 \Psi^*(\eg) = \Theta^*(\rho^2 id_2^*(\eg)) \in C^{\newreg,\newholdreg}(\overline{M},\A).
\]
\end{proof}

I would like to thank Olivier Biquard for the idea of using Whitney approximation to streamline this result.


\section{Regularity of the Defining Function} \label{defining}

This section culminates in the proof of Theorem \ref{maincptthm},
which we restate here.

\medskip
\noindent \textbf{Theorem \ref{maincptthm}.}
\emph{Let $\overline{M}$ be a compact $n$-dimensional $C^{\infty}$ manifold
with boundary, $n \geq 4$ and even.  For $\basereg \geq n$ and $0 < \holdreg < 1$, let
$\eg$ be a $C^{\basereg-1,\holdreg}$ conformally compact Einstein metric on
$M$.  Suppose that the conformal infinity of $\eg$ contains a metric
$h \in C^{\newreg,\newholdreg}(\partial M)$, where $\newreg \geq \basereg$ and $\newreg
\geq n+1$, and $0 < \newholdreg < 1$.  Then there is a $C^{\basereg, \holdreg}$ diffeomorphism
$\Psi: \overline{M} \longrightarrow \overline{M}$ which restricts to the identity on
$\partial M$ such that $\Psi^*(\eg)$ is $C^{\newreg,\newholdreg'}$ conformally compact for
some $\newholdreg'$, $0 < \newholdreg' \leq \newholdreg$.
}

\medskip
As discussed in the introduction, Theorem \ref{mainglobalthm} does not result in a $C^{\newreg, \newholdreg}$ conformally compact metric because while
$\rho^2 \eg$ is in $C^{\newreg,\newholdreg}$, the defining function that is used to generate the compactification need not be $C^{\newreg+1, \newholdreg}$.
Since the regularity for the
compactification produced in Theorem \ref{mainglobalthm} relies on a
change of smooth structure, and the defining function used wasn't
necessarily $C^{\newreg+1,\newholdreg}$ up to the boundary in the
first place, it is not immediately clear what can be said about the
defining function.

We will analyze this problem via the singular Yamabe problem, and we
will find that Theorem \ref{maincptthm} follows as a corollary to the results of our analysis.

\subsection{The Singular Yamabe Problem}
Let $\overline{M}$ be a $C^{\newreg+1,\newholdreg}$ $n$-dimensional
compact manifold with boundary, where $\newreg \geq 2$ and $0 <
\newholdreg < 1$.  Our discussion here does not depend on whether $n$
is even or odd, so until the proof of Theorem
\ref{maincptthm}, $n$ need not be even.  Let $\mathring{g}$ be a
$C^{\newreg,\newholdreg}$ conformally compact metric, with $\newreg
\geq 2$.  In the singular Yamabe problem, the goal is to find a
function $u$ such that the metric $g_+ =
u^{\frac{4}{n-2}}\mathring{g}$ has constant scalar curvature $S_+$,
which we take to be $-n(n-1)$.  In this setting, the conformal change
formula for scalar curvature becomes
\[
\mathring{\Delta}u - \frac{(n-2)}{4(n-1)}\mathring{S}u +
              \frac{(n-2)}{4(n-1)}S_+u^{\frac{n+2}{n-2}} = 0.
\]
This problem is studied as a special case of equation (7.1.1) in
\cite{AndChr}, and there it is shown that there is a function $u$ that
solves this equation, and hence the associated singular Yamabe
problem.  Moreover, the solution is unique in the class of uniformly
bounded, uniformly bounded away from zero, $C^2$ functions on $M$.
The boundary regularity of $u$ is also studied, and it is this
regularity that we will use to prove Theorem \ref{maincptthm}.

\subsection{Regularity in General}

At this point, we state the regularity result that we will be using.
This is a special case of Theorem 7.4.7 in \cite{AndChr}.  Their
result is stated for $C^{\infty}$ manifolds, but it can be checked
that their results are valid for manifolds with lower regularity
smooth structures.  Also see \cite{Andersson} for a similar result,
and \cite{Mazzeo} for a more general discussion of the singular Yamabe
problem.

\begin{thm}[Andersson, Chru\'{s}ciel] \label{Andreg}
Let $\overline{M}$ be an $n$-dimensional $C^{\newreg+1,\newholdreg}$
compact manifold with boundary, $\newreg \geq 2$, $ 0 <
\newholdreg < 1$.  Suppose $\mathring{g}$ is
$C^{\newreg,\newholdreg}$ conformally compact and 
let $u \in C^2(M)$ be a function that is uniformly
bounded, uniformly bounded above zero, and that satisfies

\[
\mathring{\Delta}u - \frac{(n-2)}{4(n-1)}\mathring{S}u 
	  - \frac{n(n-2)}{4}u^{\frac{n+2}{n-2}} = 0.
\]
Then for some $0 < \newholdreg' < 1$,
\[
u \in \left\{ \begin{array}{ll}
              C^{\newreg,\newholdreg'}(\overline{M})   & (k \leq n-1) \\
	      C^{n-1,\newholdreg'}(\overline{M}) & (k \geq n).
              \end{array} \right.
\]
Moreover, if $\newreg > n$ then there is a function $\Phi \in
C^{\newreg,\newholdreg'}(\overline{M})$ and a sequence of
functions
\[
\phi_j \in \bigcap_{i = 0}^n
y^{jn-i}C^{\newreg-n+i,\newholdreg'}(\overline{M}),\ j =
1,\ldots,N
\]
where $N$ is the smallest integer such that $N >
\newreg/n$ and $y$ is a $C^{\newreg+1,\newholdreg}$ defining function for
$\partial M$, such that
\begin{equation} \label{u_expr}
u = \Phi + \sum_{j=1}^N \phi_j \log^j(y)
\end{equation}
in $\overline{M}$.  Finally, if $y^{-n} \phi_1|_{\partial M} = 0$,
then all the $\phi_j$ can be taken to be zero and $u \in
C^{\newreg,\newholdreg'}(\overline{M})$.
\end{thm}

We note here that this theorem is global, and currently no
local analogue is known.  It is also natural to guess that the theorem
is also true setting $\newholdreg' = \newholdreg $, and that \eqref{u_expr} is
valid when $\newreg = n$.

\subsection{Regularity when $\eg$ is Einstein}

In our setting, the metric $\eg$ is Einstein, which
is stronger than simply having constant scalar curvature.  We recall
that for us, $\eg$ being Einstein means that $\eRic = -(n-1)\eg$.
The following proposition shows that in this case, the $\phi_j$ are
all zero.

\begin{prop} \label{Einstusmooth}
Let $\overline{M}$ be an $n$-dimensional $C^{\newreg+1,\newholdreg}$
compact manifold with boundary, where $\newreg \geq 2, \newreg \neq
n$, and $0 < \newholdreg < 1$, and let $\mathring{g}$ be
$C^{\newreg,\newholdreg}$ conformally compact.  Suppose $\eg =
u^{\frac{4}{n-2}} \mathring{g}$ is Einstein, where $u \in C^{2}(M)$ is
uniformly bounded and uniformly bounded above zero.  Then $u \in
C^{\newreg,\newholdreg'}(\overline{M})$ for some $\newholdreg'$, $0 <
\newholdreg' < 1$.
\end{prop}

\begin{proof}
Note that for $\newreg < n$, the Einstein condition is not necessary, and
the result is immediate from Theorem \ref{Andreg} if we only know that
$\eS = -n(n-1)$.  For $\newreg > n$, we will use the expression
\eqref{u_expr} for $u$ in Theorem \ref{Andreg} and use the fact that
$\eg$ is Einstein to show that $(y^{-n}\phi_1)|_{\partial M} = 0$.

Working in a boundary adapted coordinate system $\{x^{\alpha}\}$ with
$x^0 = y$, we let $u^{\frac{2}{n-2}} = e^{v}$ so that $\eg =
e^{2v}\mathring{g}$.  To make use of the condition that $\eg$ is
Einstein, consider the way that the Ricci tensor changes under
conformal change of the metric:
\[
\mathring{R}ic_{\alpha \beta}
      = (\mathring{\Delta} v) \mathring{g}_{\alpha \beta}
        + (n-2)(\mathring{\nabla}_{\alpha} \mathring{\nabla}_{\beta} v
                + |dv|^2_{\mathring{g}}\mathring{g}_{\alpha \beta}
		- v_{\alpha} v_{\beta})
        + (\eRic)_{\alpha \beta}.
\]
For our purposes, it will be enough to consider the trace free part of
$\mathring{R}ic$, and we note that since $\eg$ is Einstein, the trace
free part of the last term above is zero.  We denote the trace free
part of a tensor $T$ by $\tf(T)$, so in our case we have
\[ 
\tf(\mathring{R}ic_{\alpha \beta})
      = (n-2)\,\tf(\mathring{\nabla}_{\alpha} \mathring{\nabla}_{\beta} v
                           - v_{\alpha} v_{\beta}).
\]
Our analysis will actually apply to $y\, \tf (\mathring{R}ic_{\alpha
  \beta})$, so we focus on the equation
\begin{equation} \label{ytfRic}
y\,\tf(\mathring{R}ic_{\alpha \beta})
      = (n-2)\,\tf(y\,\mathring{\nabla}_{\alpha} \mathring{\nabla}_{\beta} v
                           - y\,v_{\alpha} v_{\beta}).
\end{equation}

Focusing first on the left hand side of this equation, we determine
how $\mathring{R}ic$ behaves at the boundary by making another
conformal change, namely $\csg = y^2
\mathring{g}$.  We have
\[
\mathring{R}ic_{\alpha \beta} 
       = \csRic_{\alpha \beta}
         + y^{-1}\Bigl((n-2)\cscoder_{\alpha}\cscoder_{\beta} y
                  +\csLap y\, \csg_{\alpha\beta}\Bigr)
         + y^{-2}\Bigl((1-n)|dy|_{\csg}^2 \csg_{\alpha\beta}\Bigr).
\]
Hence
\[
y\,\tf(\mathring{R}ic_{\alpha\beta})
      = \tf(y \csRic_{\alpha\beta}
            + (n-2)\cscoder_{\alpha}\cscoder_{\beta} y)
\]
which is $C^{\newreg-2,\newholdreg}(\overline{M})$, since $\csg$ is in
$C^{\newreg,\newholdreg}(\overline{M})$.  Hence $y\,\tf(\mathring{R}ic_{\alpha
\beta}) \in C^{n-1,\newholdreg}(\overline{M})$, since $\newreg > n$.

For the right hand side of \eqref{ytfRic}, we use \eqref{u_expr} to
better understand the derivatives of $v$.  We start by making a few
reductions with regard to \eqref{u_expr}.  First, by making a
conformal change to absorb $\Phi$ into $\mathring{g}$, we may say
$\Phi = 1$.  Also, we have $\phi_1 = (n-2)y^n \psi/2$, where $\psi \in
\bigcap_{i = 0}^n
y^{-i}C^{\newreg-n+i,\newholdreg'}(\overline{M})$.  Finally, let
$w= \sum_{j=1}^N \phi_j \log^j(y)$.  We can then write $v$ as follows:
\begin{eqnarray*}
v & = & \frac{2}{n-2}(w + \log(1+w) - w) \\
  & = & y^n \log(y)\psi
        + \frac{2}{n-2}\left(
          \sum_{j=2}^N \phi_j \log^j(y) + \log(1+w) - w \right)\\
  & = & y^n \log(y)\psi + f,
\end{eqnarray*}
where $f = \frac{2}{n-2}\left(\sum_{j=2}^N \phi_j \log^j(y) +
\log(1+w) - w\right)$.  By Theorem \ref{Andreg}, we know that $u \in
C^{n-1}(\overline{M})$ and $\phi_j \in \bigcap_{i = 0}^n
y^{jn-i}C^{\newreg-n+i,\newholdreg'}(\overline{M})$.  From these facts
we find that $f \in o(y^n) \cap C^{n-1,\newholdreg'}(\overline{M})$.
Note that this implies $\partial^l f \in o(y^{n-l})$ for $l \leq n-1$.

Differentiating and using the characterization of $\psi$ above,
we have
\begin{equation} \label{va}
v,_{\alpha} = n y^{n-1} y,_{\alpha}\,\log(y)\psi + O(y^{n-1})
\end{equation}
and
\begin{equation} \label{vab}
v,_{\alpha \beta} = n(n-1) y^{n-2} y,_{\alpha}\,y,_{\beta}\,
                       \log(y)\psi + O(y^{n-2}).
\end{equation}
Our next step is to write out the second covariant derivative of $v$
with respect to $\mathring{g}$ in terms of $\csg = y^2 \mathring{g}$.
In doing so, we use the following transformation rule for the
Christoffel symbol:
\[
\mathring{\Gamma}_{\alpha \beta}^{\eta}
       = \csChr_{\alpha \beta}^{\eta}
               - y^{-1}(y,_{\alpha} \delta_{\beta}^{\eta}
                        + y,_{\beta} \delta_{\alpha}^{\eta}
                        - y,_{\mu} \csg^{\eta \mu} \csg_{\alpha \beta}).
\]
With this, \eqref{va}, and \eqref{vab}, we have
\begin{eqnarray*}
\mathring{\nabla}_{\alpha} \mathring{\nabla}_{\beta} v
   & = & (v,_{\alpha \beta}
         - v,_{\eta}\mathring{\Gamma}_{\alpha \beta}^{\eta}) \\
   & = & y^{n-2} \log(y) \psi \Bigl((n^2+n)y,_{\alpha} y,_{\beta}
                            -  |dy|_{\csg}^2 \csg_{\alpha \beta}\Bigr)
+ O(y^{n-2}).
\end{eqnarray*}
Since we will be inserting this into \eqref{ytfRic}, we multiply this
equation by $y$, and observe that the trace free part of the second
term is zero.  Also, we write $y,_{\alpha} y,_{\beta} =
\delta^0_{\alpha} \delta^0_{\beta}$, and note that
$\tf(\delta^0_{\alpha} \delta^0_{\beta})$ is not zero because as a map
from $T\overline{M}$ to $T^*\overline{M}$, the transformation
$\delta^0_{\alpha} \delta^0_{\beta}$ is rank one and so cannot be a
multiple of the metric.  Finally, $v_{\alpha}v_{\beta} \in
O(y^{n-2})$, so \eqref{ytfRic} reduces to
\[    
y\,\tf(\mathring{R}ic_{\alpha \beta})
   = (n-2)(n^2+n)y^{n-1}\psi \log(y)\,
         \tf(\delta_{\alpha}^{0}\delta_{\beta}^{0}) + O(y^{n-1}).
\]
Combining this with our analysis of the left hand side of
\eqref{ytfRic}, we conclude
\[
(n-2)(n^2+n)y^{n-1}\psi \log(y)\,
         \tf(\delta_{\alpha}^{0}\delta_{\beta}^{0}) + O(y^{n-1})
   \in C^{n-1,\newholdreg}(\overline{M}).
\]
This requires that $\psi = 0$ at every point in $\partial M$, and so
by Theorem \ref{Andreg}, $u \in C^{\newreg,\newholdreg'}(\overline{M})$
\end{proof}

\subsection{Regularity of the Defining Function}

Using Proposition \ref{Einstusmooth}, we have the following proposition.

\begin{prop} \label{Einstdefsmooth}
Let $\overline{M}$ be an $n$-dimensional $C^{\newreg+1,\newholdreg}$
compact manifold with boundary, $\newreg \neq n$, with
$C^{\newreg+1,\newholdreg}$ defining function $y$.  Suppose $\eg$ is
an Einstein metric on $M$ with the property that $\csg = \csdeff^2
\eg$ extends to a $C^{\newreg,\newholdreg}$ metric on $\overline{M}$,
where $\csdeff \in C^{2}(M) \cap C^1(\overline{M})$ is a defining
function for $\partial M$.  Then $\rho/y \in
C^{\newreg,\newholdreg'}(\overline{M})$, for some $\newholdreg'$, $0 <
\newholdreg' < \newholdreg$, and therefore $\eg$ is
$C^{\newreg,\newholdreg'}$ conformally compact.
\end{prop}

\begin{proof}
We need to show that $y^2 \eg$ extends to a $C^{\newreg,\newholdreg'}$
metric on $\overline{M}$.  We will use Proposition \ref{Einstusmooth}
and some manipulation of the defining functions to achieve the result.
Observe that $\eg = u^{\frac{4}{n-2}}\mathring{g}$, where $u =
(\rho/y)^{\frac{2-n}{2}}$ and that $\mathring{g} = y^{-2}\csg$ is
$C^{\newreg,\newholdreg}$ conformally compact.  Moreover, $u$
satisfies all the hypotheses in Proposition \ref{Einstusmooth}, and so
$\rho/y \in C^{\newreg,\newholdreg'}(\overline{M})$.  Now we observe
that
\[
y^2 \eg = y^2 \csdeff^{-2} \csdeff^2 \eg = u^\frac{4}{n-2} \csg,
\]
and both $u$ and $\csg$ are $C^{\newreg,\newholdreg'}$ on $\overline{M}$.
\end{proof}

Now we prove Theorem \ref{maincptthm} as a corollary.

\begin{proof}[Proof of Theorem \ref{maincptthm}]
The proof of Theorem \ref{mainglobalthm} provides us with a $C^{\newreg+1,\newholdreg}$ atlas $\B$ which is
$C^{\basereg, \holdreg}$ related to $\A$, a $C^{\newreg+1,\newholdreg}$ diffeomorphism
$\Theta: (\overline{M}, \A) \longrightarrow (\overline{M},\B)$, and a defining function
$\csdeff \in C^{\basereg-1,\holdreg}(\overline{M}, \B)$.  In this setting, by Proposition \ref{Einstdefsmooth}, there is a defining function $y \in  C^{\newreg+1, \newholdreg}(\overline{M}, \B)$ such that
$y^2 id_2^*(\eg) \in C^{\newreg, \newholdreg'}(\overline{M}, \B)$.  Similarly to the proof of Theorem \ref{mainglobalthm}, we find that $\Psi^*(\eg)$ is $C^{\newreg, \newholdreg'}$ conformally compact.  Indeed, $\Theta^*(y) \in C^{\newreg+1, \newholdreg}(\overline{M}, \A)$ and
\[
\Theta^*(y)^2 \Psi^*(\eg) = \Theta^*(y^2 id_2^*(\eg)),
\]
which is in $C^{\newreg,\newholdreg'}(\overline{M},\A)$.
\end{proof}


\section{An Alternative Approach:  Viewing the Boundary Value Problem
  as a System} \label{syschap}

The boundary value problem \eqref{system}--\eqref{order_2l} can
naturally be interperated as an elliptic system with accompanying
boundary conditions.  Here, we introduce the framework for general
elliptic boundary value problems following \cite{ADN2} and we show
that, by incorporating \eqref{difA}, a mild adjustment of \eqref{system}--\eqref{order_2l} fits this framework.  We finish
with a discussion of the applicability of various regularity results
for such problems, as provided in \cite{ADN2} and \cite{Morrey}.

\subsection{General Elliptic Boundary Value Problems}

Here, we discuss the framework for general elliptic
boundary value problems, following closely the treatment in
\cite{ADN2}.  We focus on linear systems, but we note that Morrey
treats the nonlinear case in \cite{Morrey}.  We also restrict our
attention to boundary regularity, although \cite{ADN2} and
\cite{Morrey} deal with interior regularity as an essential precursor
to their boundary regularity results.

Throughout this discussion, Let $U$ be a coordinate domain of a $C^{\infty}$
manifold with boundary $\overline{M}$, and let $D = U \cap \partial M$
be the boundary portion of $U$.  Let $\{x^{\alpha}\}$ be boundary
adapted coordinates on $U$.

On $U$, we consider the system
\begin{equation} \label{gensys}
L_{st}(\partial) u^t = F_s,
\end{equation}
where $L_{st}(\partial)$ are the components of an $N \times N$ matrix of
differential operators, $u^t$ are the components of a vector of
unknowns, and $F_s$ are inhomogeneous terms.

On $D$, we consider
\begin{equation} \label{genbdy}
B_{rt}(\partial) u^t = \phi_r,
\end{equation}
where $B_{rt}(\partial)$ are the components of an $M \times N$ matrix of
differential operators and $\phi_q$ are inhomogeneous terms.  Note
that $N$ is determined according to how many unknowns there are, while
for now, $M$ is not yet determined.

\subsubsection{Weights and Ellipticity}

In order to determine whether or not our system is elliptic, and
whether or not our boundary data are appropriate, we need to know what
part of each differential operator should be counted as its principal
part.  To this end, we introduce integer weights that we attach to
each function $u^t$, and to each equation in the system and the
boundary conditions.

We label the weights for the functions $u^t$ by $w(u^t)$, and for the
rows of the system and boundary matrices by $w(L_s)$ and $w(B_r)$
respectively.  Letting $\ord(\mathcal{O})$ be the order of the
operator $\mathcal{O}$, the goal is to find values for these weights
such that $\ord(L_{st}(\partial)) \leq w(L_s)+w(u^t)$ and such that
$\ord(B_{rt}(\partial)) \leq w(B_r)+w(u^t)$.  Given a collection of
weights that satisfy these conditions, we can add any integer to all
the $w(u^t)$, while subtracting the same integer from $w(L_s)$ and
$w(B_r)$ to find another solution.  Hence, to eliminate this freedom,
we requre that $w(L_s) \leq 0$, with the largest such weight equal to
zero.  While other choices could be made in this regard, this choice
follows the convention used in both \cite{ADN2} and \cite{Morrey}.

With these weights in hand, we define the \emph{principal part} of
$L_{st}(\partial)$ or $B_{rt}(\partial)$ to be the term of order
exactly $w(L_s)+w(u^t)$ or $w(B_r)+w(u^t)$ respectively.
We denote this by $L_{st}'(\partial)$, respectively
$B_{rt}'(\partial)$.  If a component of either matrix has no term of
the given order, then its principal part is 0.  We also denote the
\emph{principal symbol} by $L_{st}'(\xi)$, respectively
$B_{rt}'(\xi)$, where $\xi$ is a (real) covector replacing ``$\partial$''.

We say the system $L$ is \emph{elliptic} if
$\det(L'(\xi)) \neq 0$ for any (real) non-zero covector.  Letting $m =
\frac{1}{2}\deg(\det L'(\xi))$, we say $L$ is \emph{uniformly
  elliptic} if there is a positive constant $a$ such that
$a^{-1}|\xi|^{2m} \leq |\det(L'(\xi))| \leq a |\xi|^{2m}$.

\subsubsection{Boundary Equations and the Complementing Condition}

The boundary conditions in a given boundary value problem need to be
appropriate for the system in order to make the problem well posed.
There are two conditions that must be satisfied in order for this to
happen.  First, with regard to the size of $B$, the number of
conditions $M$ must be equal to $m$ as defined above.  Second,
the boundary conditions must satisfy a
``complementing condition'' depending on an algebraic relationship
between the principal symbols of $L(\partial)$ and $B(\partial)$.

At each point $p \in D$, consider the characteristic equation
$\det(L'(\xi + \tau \nu)) = 0$ where $\xi \in T_p^*(\partial M)$ is
nonzero and $\nu$ is the inward pointing unit normal covector.  If $L$
is elliptic, then as a polynomial in the complex variable $\tau$,
this characteristic equation will have $m$ roots $\tau_r$, $1 \leq r
\leq m$, with positive imaginary part.  Consider the polynomial
\[
M^+(\tau) = \prod_{r=1}^{m} (\tau - \tau_r).
\]
We say that \emph{$B(\partial)$
satisfies the complementing condition for $L(\partial)$} if for each
$\xi$, the rows
of $B'(\xi+\tau \nu) \adj L'(\xi +\tau \nu)$ are linearly independent
modulo $M^+$, where $\adj L'$ is the matrix
adjoint to $L'$ (not the conjugate transpose).  That is to say, if
\[
c^r B_{rt}'(\xi+\tau \nu) (\adj L')^{tq}(\xi +\tau \nu)
           = P^q(\tau) M^+(\tau),
\]
where $c^r \in \mathbb{C}$, and $P^q$ are polynomials in $\tau$, then
in fact all the $c^r$ are zero.

Occasionally, in order to check the complementing condition, we can simplify our task
by observing some general facts about linear independence when working
with polynomials:

\begin{lemma} \label{linindsimp}
Let $\{v_i\}$ be a set of vectors with components that are polynomials
in one complex variable, and let $P$ be a polynomial in one complex
variable.  Let $\tau_0$ be a common root for all the components of all
$v_i$.  Then, one of two situations will occur:
\begin{enumerate}
\item If $\tau_0$ is a root of $P$, then $\{v_i\}$ is linearly
  independent modulo $P$ if and only if $\{v_i/(\tau-\tau_0)\}$ is linearly
  independent modulo $P/(\tau-\tau_0)$.
\item If  $\tau_0$ is not a root of $P$, then $\{v_i\}$ is linearly
  independent modulo $P$ if and only if $\{v_i/(\tau-\tau_0)\}$ is linearly
  independent modulo $P$.
\end{enumerate}
\end{lemma}

The proof of this lemma is straightforward.

\subsection{Application to the Current Problem}

Our goal now is to show that the boundary value problem
\eqref{system}--\eqref{order_2l} is indeed an elliptic system with boundary equations that satisfy the
complementing condition.  We will find that the system is uniformly elliptic.  The boundary equations
need to be altered however, since \eqref{difA} must be incorporated in order that
\eqref{order_0} -- \eqref{order_2l} satisfy the complementing condition.

\subsubsection{The System and Uniform Ellipticity}

Before looking at the boundary equations, we can show that \eqref{system} is uniformly elliptic.
Our first step is to do some relabeling.  The unknowns in our system
are the components of our constant scalar compactification $\csg$,
so by symmetry we have $N = n(n+1)/2$ unknowns, and for what follows,
we choose the upper triangular components as our representatives so
that we may order them in a well defined way.  To do so, we use a
non-standard lexicographic ordering for the components $\csg_{\alpha
\beta}$.  It is non-standard in that we take ``0'' to be larger than
other integers so that the components of the metric of the form
$\csg_{0 \alpha}$ come last in the list, with $\csg_{00}$ ending the
sequence.  We may then denote the functions $\csg_{\alpha \beta}$ by
$u^t$, $1 \leq t \leq N$.

Focusing on \eqref{system}, we see that our system is quasi-linear,
but for our purposes we may take it to be linear since we already
have a solution.  With this in mind, we rewrite \eqref{system} as
follows:
\[
\left\{\begin{array}{rcl}
L_{st}u^t & = & \delta_{st} \D_{\frac{n}{2}} u^t \\
F_s & = & \P(\csinv, \partial^{n-1} \csg)
\end{array}\right.
\]
for $1 \leq s \leq N$.

We now have the following:

\begin{prop}
Interpreted as above, \eqref{system} is uniformly elliptic.
\end{prop}

\begin{proof}
For weights we take $w(u^t) = n$ and
$w(L_s) = 0$, and so $L_{st}'(\partial) =  L_{st}(\partial)$.  The principal symbol for $\D$ is $|\xi|_{\csg}^2$ and more generally, the principal symbol for $\D_l$ is $|\xi|_{\csg}^{2l}$.  Hence,
for our system, we have
\[
\det(L'(\xi)) = |\xi|_{\csg}^{nN},
\]
from which we may conclude that our system is uniformly elliptic.
\end{proof}

\subsubsection{Boundary Equations and the Complementing Condition}

Focusing on the boundary equations, our first task is to incorporate \eqref{difA} into \eqref{order_0} -- \eqref{order_2l}.  We cannot simply add
\eqref{difA} to the list of boundary equations since the number of equations we need is determined by the degree of the system.  Instead, we will modify the equations of order 1 by differentiating them and substituting certain terms using \eqref{difA}.  These new equations together with the rest will be shown to satisfy the complementing condition.

To incorporate \eqref{difA} into the equations of order 1, we calculate the normal derivative
of \eqref{order_1} to get
\begin{equation*}
\csg^{\gamma 0} \csg^{\eta \beta} \csg_{\alpha \beta},_{\eta \gamma}
	- \frac{1}{2}\csg^{\gamma 0} \csg^{\eta \beta} \csg_{\eta \beta},_{\alpha \gamma}
	= \P(\csinv, \partial \csg).
\end{equation*}
When $\alpha = 0$, this becomes
\begin{equation}\label{Norder_1/0}
\csg^{\gamma 0} \csg^{\eta \beta} \csg_{0 \beta},_{\eta \gamma}
	- \frac{1}{2}\csg^{\gamma 0} \csg^{\eta \beta} \csg_{\eta \beta},_{0 \gamma}
	= \P(\csinv, \partial \csg).
\end{equation}
For $\alpha = i$, we can expand some of the sums and rearrange to produce
\begin{eqnarray} \label{Norder_1'}
\lefteqn{(\csg^{00})^2 \csg_{0 i},_{00} + 2 \csg^{00} \csg^{0k} \csg_{0i},_{0k}
	     + \csg^{0j} \csg^{0k} \csg_{0i},_{jk} - \csg^{0 \gamma} \csg^{0j} \csg_{0j},_{\gamma i} }
	     \hspace{3in} \nonumber\\
\lefteqn{+ \csg^{0 \gamma} \csg^{\eta j} \csg_{ij},_{\eta \gamma}
	      - \frac{1}{2} \csg^{0 \gamma} \csg^{jk} \csg_{jk},_{\gamma i}
	      - \frac{1}{2} \csg^{00} \csg^{0 \gamma} \csg_{00},_{\gamma i} }
	      \hspace{2.5in} \\
& = & \P(\csinv, \partial \csg) \nonumber
\end{eqnarray}
We will use \eqref{difA} to replace the first two terms in this equation.  To do this, expand the derivative of the second fundamental form in \eqref{difA} in terms of the metric and trace to get
\[
\csg^{ij} \csg_{0i},_{jk} - \frac{1}{2} \csg^{ij} \csg_{ij},_{0k}
	= \frac{(n-1)(\csg^{00})^{-1}}
                         {2(2-n)} \D \csg_{0k} 
                       + \P(\csinv, \csroot, \partial \csg, \partial_t^2 h).
\]
Expanding $\D \csg_{0k}$, rearranging, multiplying by $\csg^{00}$, and changing indices leads to
\begin{eqnarray*}
\lefteqn{(\csg^{00})^2 \csg_{0 i},_{00} + 2 \csg^{00} \csg^{0k} \csg_{0i},_{0k}} \hspace{1in} \\
	& = & \frac{2 (2-n) (\csg^{00})^2}{n-1} \csg^{jk} \csg_{0k},_{ij}
	           - \csg^{00} \csg^{jk} \csg_{0i},_{jk}
		  - \frac{(2-n) (\csg^{00})^2}{n-1} \csg^{jk} \csg_{jk},_{0i} \\
	&    & \hspace{.5in} +  \P(\csinv, \csroot, \partial \csg, \partial_t^2 h).
\end{eqnarray*}
Using this to replace the first two terms on the left hand side of \eqref{Norder_1'}, we have
\begin{eqnarray} \label{Norder_1',w/difA}
\lefteqn{\frac{2 (2-n) (\csg^{00})^2}{n-1} \csg^{jk} \csg_{0k},_{ij}
	   - \csg^{00} \csg^{jk} \csg_{0i},_{jk}
	   - \frac{(2-n) (\csg^{00})^2}{n-1} \csg^{jk} \csg_{jk},_{0i}
	   + \csg^{0j} \csg^{0k} \csg_{0i},_{jk}} \hspace{3in} \nonumber \\
\lefteqn{- \csg^{0 \gamma} \csg^{0j} \csg_{0j},_{\gamma i}
	   + \csg^{0 \gamma} \csg^{\eta j} \csg_{ij},_{\eta \gamma}
	   - \frac{1}{2} \csg^{0 \gamma} \csg^{jk} \csg_{jk},_{\gamma i}
	   - \frac{1}{2} \csg^{0 \gamma} \csg^{00} \csg_{00},_{\gamma i}} \hspace{2.5in} \\
& = & \P(\csinv, \csroot, \partial \csg, \partial_t^2 h), \nonumber
\end{eqnarray}

which, with \eqref{Norder_1/0}, are used in place of the equations of order 1 to give us the following:

\begin{prop} \label{bdycomp}
The equations \eqref{order_0}, \eqref{Norder_1/0}, \eqref{Norder_1',w/difA},
\eqref{order_2} -- \eqref{order_2l} satisfy the complementing condition for the system \eqref{system}.
\end{prop}

After introducing some organizing notation and calculating the polynomial $M^+$ and the matrix
adjoint to $L_{st}'(\xi + \tau \nu)$, we prove a technical lemma to help us analyze these equations.  Finally, we prove this proposition using an inductive argument.

We have a total of $M = nN/2$ boundary equations so for the components
$B_{rt}(\partial)$ we let $1 \leq r \leq M$ and as above, we
have $1 \leq t \leq N$.  These
conditions are grouped by the order of the equations and we will call the set
of rows in $B_{rt}(\partial)$ corresponding to the equations of order
$s$ the \emph{block of order $s$}.  Additionally, we will refer to the rows corresponding to
\eqref{order_2} as the block of order 2(a), while the rows corresponding to \eqref{Norder_1/0} and
\eqref{Norder_1',w/difA} will be denoted the block of order 2(b).  As with the system, we relegate
all but the highest order parts of each equation to the inhomogeneous term.
With $w(u^t) = n$ as indicated above, we let $w(B_r) = s-n$
for any row $B_r(\partial)$ in a block of order $s$.  Hence we have $B'(\partial) = B(\partial)$.  

Working out how each of the boundary equations are represented, we find that $B(\partial)$ has the following properties:  The leftmost $(N-n) \times (N-n)$ minor in the block of order 0 is the identity matrix $I_{N-n}$.  The remaining entries in the block of order 0 are equal to zero.  The leftmost $(N-n) \times (N-n)$ minor in the block of order 2(a) is $\D I_{N-n}$.  The remaining entries in the block of order 2(a) are equal to zero.  The blocks of order $2l$, $l \geq 2$ are of the form $\D_l I_N$.  The rightmost $n \times n$ minor in the block of order 3 is lower triangular and, with the exception of
the last row, it is diagonal.  The operator on the diagonal is $\csg^{\eta 0} \partial_{\eta} \D$, while the operator in the $j$th entry of the last row is $\csg^{\eta j} \partial_{\eta} \D$.  We will find that the remainder of this block will not enter into our analysis.  Finally, the block of order 2(b) proves to be complicated.  Like the block of order 3, we will find that all we need to study is the rightmost $n \times n$ minor.  This will be studied more closely later.

We need the polynomial $M^+$ and the matrix
adjoint to $L_{st}'(\xi + \tau \nu)$ in order to check the complementing condition.
In our case, we have $\det(L'(\xi + \tau \nu)) = |\xi + \tau \nu|_{\csg}^{nN}$.
This function has two complex roots, $\tau^+$ and $\tau^- =
\overline{\tau^+}$, each with multiplicity $nN/2$.
Taking $\tau^+$ to be the root with positive imaginary part
we have
\[
M^+(\tau) = (\tau - \tau^+)^{\frac{n}{2}N}.
\]
Also, calculating the matrix adjoint we have
\[
\adj (L')^{tq}(\xi + \tau \nu) = |\xi +\tau \nu|_{\csg}^{n(N-1)}
                                 \delta^{tq}.
\]
Note that, as with the determinant of $L'(\xi+\tau \nu)$, each
component here has the two roots $\tau^+$, and $\tau^-$, but this
time their multiplicities are each $n(N-1)/2$.

Now, through Lemma \ref{linindsimp}, the complementing condition simplifies
in that all of the terms in $\adj L'(\xi +\tau \nu)$ may be reduced to
$\delta^{tq}$ and in doing so, $M^+$ is reduced to $(\tau -
\tau^+)^{\frac{n}{2}}$.  Hence, the only effect that $\adj L'$ has is
to raise an index on $B'$ and to keep notation simple, we lower it
back down.  Therefore, verifying the complementing condition amounts to checking that the rows of
$(B')(\xi+\tau\nu)$ are linearly independent modulo
$(\tau-\tau^+)^{\frac{n}{2}}$.  This follows from a careful analysis of
the interaction among the various rows in $B'(\xi + \tau \nu)$.  For
the most part, this analysis is straightforward.  The only real
complication is that the rows composing the block of order 2(b) do not
have a simple structure.  We will find that the rightmost $n \times n$ minor in this block, which we label
$\F$, is all we need to analyze.  To help in this regard, we have the following lemma:

\begin{lemma} \label{order_1_block}
In a neighborhood of a point $p \in D$ where $\csg = id$, the matrix $\F$ is invertible
when $\tau = \tau^+$.
\end{lemma}

We note that by making a linear change to the coordinates in Proposition \ref{hcoordsprop}, we can guarantee that the metric is equal to the identity at a point of our choice, so this condition does not have a significant impact on the applicability of these results.

\begin{proof}
It is sufficient to analyze $\F$ when $\csg_{\alpha \beta} = \delta_{\alpha \beta}$.  In this case, when
$\tau = \tau^+$, the matrix simplifies to 
\[
\left(\begin{array}{cc}
A & v \\
w & b
\end{array} \right),
\]
where $A$ is an $(n-1) \times (n-1)$ matrix with $A_{ij} = -\delta_{ij}|\xi|^2+\frac{2(2-n)}{n-1}\xi_i \xi_j$,
$v$ is a column vector with $v_i = -\frac{1}{2} \xi_i \tau^+$, $w$ is a row vector with $w_j = \xi_j \tau^+$, and $b = \frac{1}{2}(\tau^+)^2$.

Performing column operations on this matrix, we get a new matrix
\[
\left(\begin{array}{cc}
A' & v' \\
w' & b'
\end{array} \right),
\]
where $A'_{ij} = -\delta_{ij}|\xi|^2+\frac{3-n}{n-1}\xi_i \xi_j$, $v'_i = -\frac{1}{2} \xi_i$, $w'_j = 0$ and
$b' = \frac{1}{2}\tau^+$.

Since all but the last entry in the bottom row are zero, invertibility of the entire matrix depends on invertibility of $A'$.  Note that $A' = \frac{3-n}{n-1} \xi \xi^T - |\xi|^2 I_{(n-1)}$, so $A'$ is invertible as long as $|\xi|^2$ is not an eigenvalue of $\frac{3-n}{n-1} \xi \xi^T$.  It is easy to verify that the eigenvalues of $\frac{3-n}{n-1} \xi \xi^T$ are $\frac{3-n}{n-1} |\xi|^2$ (with multiplicity 1) and 0 (with multiplicity $n-2$).  Therefore, $A'$ is invertible, and so the original matrix $\F$ is invertible as well.
\end{proof}

Using this lemma to help deal with the block of order 2(b), we have the following:

\begin{proof}[Proof of Proposition \ref{bdycomp}]
As mentioned above, we need to show that the rows of $B'(\xi+\tau\nu)$ are linearly independent
modulo $(\tau - \tau^+)^{\frac{n}{2}}$.  That is, given the following:
\begin{equation} \label{lin_combo}
c^r (B')_{rt}(\xi+\tau\nu) = P_t(\tau)(\tau-\tau^+)^{\frac{n}{2}},
\end{equation}
where $c^r \in \mathbb{C}$ and $P_t$ are polynomials in one complex
variable, then in fact all $c^r$ are zero.  Letting $C^s$ be the (row) vector of coefficients $c^r$ corresponding to the rows in the block of order $s$, we will proceed inductively.  First we will show that $C^0$ and $C^{2(b)}$ must be zero.  Then we will show that $C^{2(a)}$ and $C^3$ must be zero.  Finally, we will show that the remaining coefficients are zero.

We note again that $\D(\xi+\tau\nu) = |\xi+\tau\nu|_{\csg}^2 = (\tau-\tau^+)(\tau-\tau^-)$.  More generally, we have $\D_l(\xi+\tau\nu) = (\tau-\tau^+)^l(\tau-\tau^-)^l$.  Also, the operators in the rightmost
$n \times n$ minor in the block of order 3 each include an $\D$ and a sum of first order operators.  Hence each of these terms take the form
$\csg^{\alpha \eta}(\tau-\tau^+)(\tau-\tau^-)(\xi_{\eta}+\tau \nu_{\eta})$ which has one real root since it is cubic and all its coefficients are real.  We can then write this minor as
$(\tau-\tau^+)(\tau-\tau^-)\T$ where $\T$ is a lower triangular $n \times n$ matrix with
$\csg^{0 \eta}(\xi_{\eta}+\tau \nu_{\eta})$ for entries on the diagonal and
$\csg^{i \eta}(\xi_{\eta}+\tau \nu_{\eta})$ for the $i$th entry of the last row.  With this we can write down \eqref{lin_combo} as follows:

\begin{eqnarray} \label{expanded_lin_combo}
\lefteqn{
(C^0\ C^{2(b)})\left(\begin{array}{cc}
			I_{N-n} & 0 \\
			*	    & \F
			\end{array}\right)
+ (\tau - \tau^+)(\tau - \tau^-)(C^{2(a)}\ C^3)\left(\begin{array}{cc}
								I_{N-n} & 0 \\
								*	    & \T
								\end{array}\right)} \hspace{2in} \nonumber \\
\lefteqn{+ C^4(\tau - \tau^+)^2(\tau - \tau^-)^2 I_N
+ . . . + C^{n-2} (\tau - \tau^+)^{\frac{n}{2} - 1}(\tau - \tau^-)^{\frac{n}{2} - 1} I_N} \hspace{2in} \\
& = & (P_1(\tau)\ \ldots \ P_N(\tau))(\tau - \tau^+)^{\frac{n}{2}} \nonumber
\end{eqnarray}

For the first step, \eqref{expanded_lin_combo} must be true in particular at $\tau = \tau^+$.  With the exception of the first term, every term is 0 at $\tau = \tau^+$.  Since the matrix in the first term is invertible at $\tau^+$ (with the help of Lemma \ref{order_1_block}), $C^0$ and $C^{2(b)}$ must be zero.

The analysis for the next term in \eqref{expanded_lin_combo} is similar.
We divide by $(\tau-\tau^+)$
and since the equation must still be true at $\tau = \tau^+$, we may
conclude that $C^{2(a)}$ and $C^3$ must be zero.  Continuing by dividing by
$(\tau-\tau^+)$ and then evaluating at $\tau = \tau^+$, we see that all
the remaining $C^s$  must be zero as well.

Therefore all $c^r$ must be zero and the rows of $B'(\xi+\tau\nu)$ are
linearly independent modulo $(\tau - \tau^+)^{\frac{n}{2}}$.
\end{proof}

\subsection{Regularity}

Once we know that our system satisfies the complementing condition, we can apply the boundary regularity results of \cite{ADN2} or \cite{Morrey} with the goal of proving Theorems \ref{mainlocalthm} and \ref{mainglobalthm}.  Unfortunately it seems that these regularity results do not give us quite the results we want.  In particular, the results in \cite{ADN2} are local, but, in our context, they also require that $\csg \in C^{n,\holdreg}(U)$.  On the other hand, there are weak regularity results in \cite{Morrey}, but they are global on a domain in $\mathbb{R}^n$ and require strong regularity transverse to the boundary.  The global nature makes it unclear how to account for the use of harmonic coordinates in this setting.  Local/global issues aside, in order to satisfy the regularity conditions at the boundary, $\csg$ is required to be in $C^{4,\holdreg}(\overline{M})$ in dimension 4 and $C^{n-1,\holdreg}(\overline{M})$ in higher dimensions, the difference arising because of the fact that in dimension 4, there is boundary data up to order 3, while in dimension $n$ greater than 4, there is boundary data only up to order $n-2$.  Therefore, neither approach improves the bootstrap approach presented here.

\bibliographystyle{amsalpha}
\bibliography{conformal}

\end{document}